\documentclass[12pt,a4paper]{article}
\usepackage{amsmath,amssymb}
\newtheorem{theorem}{Theorem}
\newtheorem{lemma}{Lemma}
\newtheorem{example}{Example}
\newtheorem{remark}{Remark}
\newcommand{\half}{\frac{1}{2}}
\newcommand{\nn}{\nonumber \\}
\def\eqn#1{(\ref{#1})}
\newcommand{\cinf}{C^\infty}       
\newcommand{\row}[3]{{#1}_{#2},\dots,{#1}_{#3}}   
\newcommand{\ii}{{\,{\rm i}\,}}
\newcommand{\ca}{{\mathcal A}}
\newcommand{\cb}{{\mathcal B}}
\newcommand{\cc}{{\mathcal C}}

\newcommand{\ch}{{\mathcal H}}
\newcommand{\ci}{{\mathcal I}}

\newcommand{\ck}{{\mathcal K}}
\newcommand{\cl}{{\mathcal L}}
\newcommand{\cm}{{\mathcal M}}

\newcommand{\cs}{{\mathcal S}}

\newcommand{\IC}{{\mathbb C}}

\newcommand{\IF}{{\mathbb F}}
\newcommand{\IH}{{\mathbb H}}
\newcommand{\II}{{\mathbb I}}

\newcommand{\IR}{{\mathbb R}}

\newcommand{\IT}{{\mathbb T}}
\newcommand{\IZ}{{\mathbb Z}}
\newcommand{\wh}{\widehat}
\newcommand{\wt}{\widetilde}
\newcommand{\pa}{\partial}
\def\bar#1{\overline{#1}}
\def\norm#1{\|#1\|}
\def\abs#1{{\vert#1\vert}}

\def\inf1{{\mathcal L}^{(1,\infty)}}

\def\cha#1{C_{#1}({\mathcal A})}
\def\ccha#1{CC_{#1}({\mathcal A})}
\def\cocha#1{C^{#1}({\mathcal A})}
\def\cochala#1{C^{#1}_{\lambda}({\mathcal A})}

\newcommand{\ot}{\otimes}
\newcommand{\op}{\oplus}

\def\ket#1{\left| #1\right\rangle}
\def\hs#1#2{\left\langle #1,#2\right\rangle}

\def\dn{\downarrow}

\def\ra{\rightarrow}

\def\lra{\longrightarrow}
\def\lla{\longleftarrow}

\numberwithin{equation}{section}
\def\lgl{\langle}
\def\rgl{\rangle}
\def\ov{\overline}
\def\ts{\times}
\def\Sup{\mathop{\rm Sup}\nolimits}
\def\fl{\forall}
\def\sbs{\subset}
\DeclareMathOperator{\im}{im}
\DeclareMathOperator{\Cliff}{Cliff}
\DeclareMathOperator{\diag}{diag}
\DeclareMathOperator{\id}{id}
\DeclareMathOperator{\Mat}{Mat}
\DeclareMathOperator{\Tr}{Tr}
\DeclareMathOperator{\Spec}{Spec}
\DeclareMathOperator{\ChC}{ch}
\DeclareMathOperator{\tr}{tr}
\DeclareMathOperator{\SU}{SU}
\DeclareMathOperator{\U}{U}
\DeclareMathOperator{\SO}{SO}
\DeclareMathOperator{\Aut}{Aut}
\DeclareMathOperator{\Hom}{Hom}
\newcommand{\Cs}{\mbox{\upshape C}\ensuremath{{}^*}}
\newcommand{\into}{\hookrightarrow}

\newcommand{\addots}{\mathinner{\mkern-1mu
\raise0.5pt\vbox{\kern7pt\hbox{.}}\mkern2mu\raise4pt\hbox{.}\mkern2mu
\raise7.5pt\hbox{.}\mkern-1mu}}
\makeindex
\begin{document}

\title{Noncommutative Spheres and Instantons}
\date{July 2003}

\author{Giovanni Landi \\[15pt]
Dipartimento di Scienze Matematiche,
           Universit\`a di Trieste\\
           Via Valerio 12/b, I-34127 Trieste, Italia \\
and INFN, Sezione di Napoli, Napoli, Italia \\
{\tt landi@univ.trieste.it}
}

\maketitle
\vspace{1 cm}

\begin{abstract}
We report on some recent work on deformation of spaces, notably deformation of spheres,
describing two classes of examples. \\  The first class of examples consists of 
noncommutative manifolds associated with the so called $\theta$-deformations  which were
introduced in \cite{cola} out of a simple analysis in terms of cycles in the $(b,B)$-complex
of cyclic homology. These examples have non-trivial global features and can be endowed with a
structure of noncommutative manifolds, in terms of a spectral triple $(\ca , \ch , D)$. In
particular, noncommutative spheres $S^{N}_{\theta}$ are isospectral deformations of usual
spherical geometries. For the corresponding spectral triple $(\cinf(S^{N}_\theta), \ch , D)$,
both the Hilbert space of spinors $\ch= L^2(S^{N},\cs)$ and the Dirac operator $D$ are the
usual ones on the commutative
$N$-dimensional sphere $S^{N}$ and only the algebra and its action on $\ch$ are deformed.
\\ The second class of examples is made of the so called quantum spheres $S^{N}_q$  which are
homogeneous spaces of quantum orthogonal and quantum unitary groups. For these spheres, there
is  a complete description of $K$-theory, in terms of  nontrivial self-adjoint idempotents
(projections) and unitaries, and of the $K$-homology, in term of nontrivial Fredholm
modules,   as well as  of the corresponding Chern characters in cyclic homology and
cohomology.
\end{abstract}

\vfill
\noindent These notes are based on invited lectures given at the {\it International Workshop
on Quantum Field Theory and Noncommutative Geometry}, November 26-30 2002, Tohoku University,
Sendai, Japan. To be published in the workshop proceedings by Springer-Verlag as Lecture
Notes in Physics.

\vfill\eject
\tableofcontents
\vfill\eject
\section{Introduction}\label{se:intro}

We shall describe two classes of deformation of spaces with particular emphasis on 
spheres.

The first class of examples are noncommutative manifolds associated with the so
called $\theta$-deformations and which are constructed naturally \cite{cola} by a simple analysis in
terms of cycles in the $(b,B)$-complex of cyclic homology. These examples have non-trivial
global features and can be endowed with a structure of noncommutative manifolds,
in terms of a spectral triple $(\ca , \ch , D)$ \cite{co94,co96}.
In particular we shall describe noncommutative spheres $S^{N}_{\theta}$ which are 
isospectral deformations of usual spherical geometries; and  we shall also show quite generally
that any compact Riemannian spin manifold whose isometry group has rank $r \geq 2$ admits
isospectral deformations to noncommutative geometries.

\bigskip

The second class of examples is made of the so called quantum orthogonal spheres $S^{N}_q$,
which have been constructed as homogeneous spaces \cite{FadResTak89} 
of quantum orthogonal groups $\SO_q(N+1)$ and quantum unitary  spheres $S^{2n+1}_q$
which are homogeneous spaces of  quantum unitary
groups $\SU_q(n+1)$ \cite{VakSoi91}.  The quantum groups $\SO_q(N+1)$ and $\SU_q(n+1)$
are $R$-matrix deformations of the usual orthogonal and unitary groups $\SO(N+1)$ and $\SU_q(n+1)$
respectively.
In fact, it has been remarked in \cite{hala} that `odd' quantum orthogonal spheres are the
same as `odd' quantum unitary ones, as it happens for undeformed spheres.

It is not yet clear if and to which extend these quantum spheres can be endowed with the structure 
of a noncommutative geometry via a spectral triple. There has been some interesting work in this
direction recently. 
In \cite{ChaPal02} a $3$-summable spectral triple was
constructed for $\SU_q(2)$; this has been thoroughly analyzed in 
\cite{co02}
in the context of the noncommutative local index formula of 
\cite{como}. 
A $2$-summable spectral triple on $\SU_q(2)$ was constructed in
\cite{ChaPal02b} together with a spectral triple on the  spheres $S^2_{qc}$
of Podle{\'s} \cite{Pod87}. Also, a `$0$-summable' spectral triple on the so called standard spheres
$S^2_{q0}$ has been given in \cite{dasi,scwa,kr}. 
Instead, on these spheres one can construct  Fredholm modules,
which provide a structure which is somewhat weaker that the one given by spectral triples. Indeed, a
Fredholm module can be though of as a noncommutative conformal structure \cite{cosute}.  This
construction for the quantum spheres $S^{N}_q$ will be described in Section~\ref{se:fmqes} closely
following the paper
\cite{hala}. 

\bigskip

All our spaces can be regard  as ``noncommutative real affine varieties". 
For such an object, $X$, the algebra $A(X)$ is a finitely presented $*$-algebra in terms of 
generators and relations.  In contrast with classical algebraic geometry, there does not in general
exist a topological point set $X$. Nevertheless, we regard $X$ as a
noncommutative space and $A(X)$ as the algebra of polynomial functions
on $X$.
In the classical case, one can consider the algebra of continuous
functions on the underlying topological space of an affine variety. If
$X$ is bounded, then this is a \Cs-algebra and is the completion of
$A(X)$. In general, one defines $\cc(X)$ to be the \Cs-algebraic
completion of the $*$-algebra $A(X)$. To construct this, one first considers the free algebra 
$F(X)$ on the same generators of the algebra $A(X)$. Then, one takes 
all possible $*$-representations $\pi$ of $F(X)$ as bounded operators on a countably
infinite-dimensional Hilbert space $\ch$. The representations are taken to be {\it admissible}, 
that is in $\cb(\ch)$ the images of the generators of $F(X)$ satisfy the same defining relations as in
$A(X)$.   For $a\in F(X)$ one  defines  
$\norm{a}  = \Sup \norm{\pi (a) } $
with $\pi$ ranging through all admissible representations of $F(X)$.  It turns out that $\norm{a}$ is
finite for $a\in F(X)$ and $\norm{\cdot}$ is a seminorm. Then $\ci :=\{a\in F(X) =0\}$ is a two-sided
ideal and one obtains a
\Cs-norm on
$F(X) / \ci$. The 
\Cs-algebra $\cc(X)$ is
the
completion of $F(X)/\ci$ with respect to this norm. The \Cs-algebra
$\cc(X)$ has the
universal property that any \hbox{$*$-morphism} from $A(X)$ to a
separable
\Cs-algebra factors through $\cc(X)$. In particular, any
\hbox{$*$-representation} of $A(X)$ extends to a representation of $\cc(X)$.

\bigskip
The word instanton in the title refers to the fact that all (in particular even) spheres come
equipped with a  projection $e\in\Mat_r(A(X)), e^2=e=e^*$, for $X=S^{N}_\theta$ and $X=S^{N}_q$.
These projection
determines the module of sections of a vector
bundle which deforms the usual monopole bundle and instanton bundle in two and four dimensions
respectively, and generalizes them in all dimensions. 

In particular on the four dimensional $S_\theta$, one can develop
Yang-Mills theory,  since there are all the
required structures, namely the algebra, the calculus and the ``vector
bundle" $e$ (naturally
endowed, in addition, with a preferred connection $\nabla$).
Among other things there is 
a basic inequality 
showing that the Yang-Mills
  action, $YM(  \nabla) ={\int \!\!\!\!\!\! -} \,  \theta ^2 \,  ds^{4} $,
(where $\theta = \nabla^2$ is the curvature, and
$ds=D^{-1}$)
  has a strictly positive
  lower bound given by the topological invariant $\varphi(e)=\int\!\!\!\!\!- \gamma (e
-\frac{1}{2}) [D, e]^4
\ ds^4$ which, for the canonical projections turns out to be just $1$: $\varphi(e)=1$. \\
In general, the projection $e$ for the spheres $S^{2n}_\theta$ satisfies self-duality 
equations
\begin{equation}
*_H e (d e)^n = \ii^n e (d e)^n ,
\end{equation}
with a suitably defined Hodge operator $*_H$ \cite{codu-vi} (see also \cite{asbo} 
and \cite{la01}).\\ 
An important problem is the construction and the classification of Yang-Mills
connections in the noncommutative situation along the line of the ADHM construction \cite{at79}.
This was done in
\cite{cori} for the noncommutative torus and in \cite{nesc} for a noncommutative $\IR^4$.

It is not yet clear if a construction of gauge theories along similar lines can be done 
for the quantum spheres $X=S^{N}_q$.

\bigskip
There has been recently an explosion of work on deformed spheres from many points of view.
The best  I can do here is to refer to \cite{da} for an overview of noncommutative and
quantum spheres in dimensions up
to four. In \cite{sit2} there is a family of noncommutative 4-spheres which satisfy the Chern
character conditions of \cite{cola} up to cohomology classes (and not just representatives). 
Additional quantum $4$-dimensional spheres together with a construction of quantum instantons
on them is in \cite{frja}. A
different class of spheres in any even dimension was proposed in \cite{bocita}.  
At this workshop T. Natsume presented an example in two dimensions
\cite{na}.

\section{Instanton Algebras}\label{inal}
In this  Section we shall describe how to obtain in a natural way noncommutative spaces 
(i.e. algebras) out of the Chern characters  of idempotents and unitaries in cyclic homology.
For this we shall give a brief overview of the needed fundamentals of the theory, following
\cite{br}. For later use we shall also describe the dual cohomological theories.

\subsection{Hochschild and cyclic homology and cohomology}\label{co-ho}
Given an algebra $\ca$, consider the chain complex 
$(C_*(\ca)=\bigoplus_n\cha{n}, b)$ with $\cha{n}=\ca^{\otimes(n+1)}$ and the boundary map $b$
defined by
\begin{eqnarray}\label{bi}
b&\::&\cha{n}\ra\cha{n-1}~, \nn
b(a_0 \ot a_1 \ot \cdots \ot a_n) &:=&
\sum_{j=0}^{n-1} (-1)^j a_0 \ot \cdots \ot a_j a_{j+1} \ot \cdots \ot a_n \nn
&&~~~~~~~~~~ + (-1)^{n} a_n a_0 \ot a_1 \ot \cdots \ot a_{n-1} \, .
\end{eqnarray}
It is easy to prove that $b^2=0$.
The {\it Hochschild homology} $HH_{*}(\ca)$ of the algebra $\ca$ is
the homology of this complex,
\begin{equation}
HH_n(\ca)  :=  H_n(C_*(\ca),b) = Z_n / B_n ~,
\end{equation}
with the cycles given by $Z_n:=\ker ( b :\cha{n}\ra\cha{n-1} )$ and the boundary given by $B_n
:=\im (b :\cha{n+1}\ra\cha{n})$.  We have another
operator which increases the degree
\begin{equation}\label{Bi}
B:\cha{n}\ra\cha{n+1}~, ~~~~~
B =  B_0 A\, ,
\end{equation}
where
\begin{eqnarray}\label{Bi0}
&& B_0 (a_0 \ot
a_1 \ot
\cdots \ot a_n ) := \II \ot a_0 \ot a_1 \ot \cdots \ot a_n
\\
&& A
(a_0 \ot a_1 \ot \cdots \ot a_n ) := \frac{1}{n + 1} \sum_{j=0}^n 
(-1)^{nj} a_j \ot a_{j+1} \ot \cdots \ot a_{j-1} \, , ~~~\label{Ai}
\end{eqnarray}
with the obvious cyclic identification $n+1 = 0$.
Again it is straightforward to check that $B^2=0$ and that $b B + B b = 0$.  

By putting together these two operators, one gets a bi-complex
$(C_*(\ca), b, B)$ with $\cha{p-q}$ in bi-degree $p,q$. 
The {\it cyclic homology} $HC_{*}(\ca)$ of the algebra $\ca$ is
the homology of the total complex $(CC(\ca), b+B)$,  
whose $n$-th term is given by $\ccha{n}  := \op_{p+q=n} \cha{p-q} = \op_{0\leq q
\leq  [n/2]} \cha{2n-q}$.
This bi-complex may be best organized in a plane diagram whose vertical arrows are associated with the
operator $b$ and whose horizontal ones are associated with the operator $B$, 

\begin{equation}\label{toco}
\begin{array}{ccccc}
  ~ & \vdots & ~    & \vdots & ~   \\
  ~ & ~      & ~    & ~      & ~   \\
b \dn & ~      & b \dn  & ~      & b \dn \\
  ~ & ~      & ~    & ~      & ~   \\
\cha{2} & ~{\buildrel B  \over \lla} & \cha{1}~ & ~{\buildrel B  \over \lla}~ & \cha{0} \\
  ~     & ~    & ~       & ~    & ~   \\
b \dn     & ~    & b \dn     & ~    & ~  \\
  ~ & ~      & ~    & ~      & ~   \\
\cha{1} & ~{\buildrel B  \over \lla}~ & \cha{0} & ~    & ~  \\
  ~     & ~    & ~       & ~    & ~   \\
b \dn     & ~    & ~       & ~    & ~ \\
  ~     & ~    & ~       & ~    & ~   \\
\cha{0} & ~    & ~       & ~    & ~ \\
  ~     & ~    & ~       & ~    & ~
\end{array} 
\end{equation}
The $n$-th term $\ccha{n}$ of the total complex is just the $n$-th (NW -- SE) diagonal in the
diagram
\eqn{toco}.   Then,
\begin{equation}
HC_n(\ca)  :=  H_n(CC(\ca), b+B) = Z^\lambda_n / B^\lambda_n~,
\end{equation}
with the {\it cyclic cycles} given by 
$Z^\lambda_n := \ker ( b+B :\ccha{n}\ra\ccha{n-1} )$ and the {\it cyclic boundaries} given by 
$B^\lambda_n :=\im (b+B :\ccha{n+1}\ra\ccha{n})$.

\begin{example} 
If $M$ is a compact manifold, the Hochschild homology of the
algebra of smooth functions $C^\infty(M)$ gives the de~Rham complex 
(Hochschild-Konstant-Rosenberg theorem),
\begin{equation}
\Omega^k_{dR}(M) \simeq HH_k(C^\infty(M))~,
\end{equation}
with $\Omega^k_{dR}(M)$ the space of de~Rham forms of order $k$ on $M$. If $d$  denotes the 
de~Rham exterior  differential, this isomorphisms is implemented by 
\begin{equation}
a_o da_1 \wedge \cdots \wedge da_k \mapsto \varepsilon_k(a_0 \ot a_1 \ot \cdots \ot da_k)
\end{equation}
where $\varepsilon_k$ is the {\it antisymmetrization map}
\begin{equation}
\varepsilon_k(a_0 \ot a_1 \ot \cdots \ot da_k) := 
\sum_{\sigma\in S_k}sign(\sigma)(a_0 \ot a_{\sigma(1)} \ot \cdots \ot da_{\sigma(k)}) 
\end{equation}
and $S_k$ is the symmetric group of order $k$. In particular one checks that $b \circ \varepsilon_k =0$. 
The de~Rham  
differential $d$ corresponds to the operator $B_*$ (the lift of
$B$ to homology) in the sense that   
\begin{equation}
\varepsilon_{k+1} \circ d = (k+1) B_* \circ \varepsilon_k ~.
\end{equation}
On the other hand, the cyclic
homology gives
\cite{co85,lo} 
\begin{equation}\label{deRa}
HC_k(C^\infty(M)) =\Omega^k_{dR}(M) / d \Omega^{k-1}_{dR}(M)\op H^{k-2}_{dR}(M) \op H^{k-4}_{dR}(M) \op
\cdots~,
\end{equation}
where $H^{j}_{dR}(M)$ is the $j$-th de~Rham
cohomology group. The last term in the sum is $H^{0}_{dR}(M)$ or $H^{1}_{dR}(M)$ according to wether $k$
is even or odd. \\
From the fact that $C^\infty(M)$ is commutative it follows that there is a natural decomposition
(the $\lambda$-decomposition) of cyclic homology in smaller pieces,
\begin{align}\label{lade}
HC_0(C^\infty(M)) &= HC_0^{(0)}(C^\infty(M))~, \nn
HC_k(C^\infty(M)) &= HC_k^{(k)}(C^\infty(M))  \cdots~ \op HC_k^{(1)}(C^\infty(M))~,
\end{align}
which is obtained by suitable idempotents $e_k^{(i)}$ which commute with the operator $B$:
$B e_k^{(i)}=e_{k+1}^{(i+1)} B$. The previous decomposition corresponds to the decomposition in
\eqn{deRa} and give a way to extract the de~Rham cohomology
\begin{equation}\label{dedeRa}
\begin{array}{ll}
HC_k^{(k)}(C^\infty(M)) = \Omega^k_{dR}(M) / d \Omega^{k-1}_{dR}(M) ~, \\ 
HC_k^{(i)}(C^\infty(M)) = H^{2i-k}_{dR}(M)~, & {\rm for} ~~~ [n/2] \leq i < n ~,  \\
HC_k^{(i)}(C^\infty(M)) = 0 ~, & {\rm for} ~~~ i < [n/2] ~,
\end{array}
\end{equation}
\end{example} 
Looking at this example, one may think of cyclic homology as a generalization of de~Rham cohomology to
the noncommutative setting.

\bigskip
A Hochschild $k$ cochain on the algebra $\ca$ is an $(n+1)$-linear functional on $\ca$ or a linear form
on $\ca^{\otimes(n+1)}$. Let $\cocha{n}=\Hom(\ca^{\otimes(n+1)},\IC)$ be the collection of such
cochains. We have a cochain complex 
$(C^*(\ca)=\bigoplus_n\cocha{n}, b)$  with a coboundary map, again
denoted with the symbol $b$, defined by
\begin{eqnarray}\label{cobi}
b&\::&\cocha{n}\ra\cocha{n+1}~, \nn
b\varphi(a_0, a_1, \cdots, a_{n+1}) &:=&
\sum_{j=0}^{n} (-1)^j \varphi(a_0, \cdots, a_j a_{j+1}, \cdots, a_{n+1}) \nn
&&~~~~~~~~~~ + (-1)^{n+1} \varphi(a_{n+1} a_0 a_1, \cdots, a_{n}) \, .~~~
\end{eqnarray}
Clearly $b^2=0$ and 
the {\it Hochschild cohomology} $HH^{*}(\ca)$ of the algebra $\ca$ is
the cohomology of this complex,
\begin{equation}
HH^n(\ca)  :=  H^n(C^*(\ca),b) = Z^n / B^n~,
\end{equation}
with the cocycles given by $Z^n:=\ker(b :\cocha{n}\ra\cocha{n+1})$ and the coboundaries given by
$B^n :=\im(b:\cocha{n-1}\ra\cocha{n})$.

A Hochschild
$0$-cocycle $\tau$ on the algebra $\ca$ is a {\it trace}, since
$\tau\in\Hom(\ca,\IC)$ and the cocycle condition is 
\begin{equation}
\tau(a_0a_1)-\tau(a_1a_0)=b\tau(a_0,a_1)=0~.
\end{equation}
The trace property is extended to higher orders by saying that an $n$-cochain $\varphi$ is  
{\it cyclic} if $\lambda \varphi = \varphi$, with
\begin{equation}
\lambda\varphi(a_0, a_1, \cdots, a_{n})=(-1)^n\varphi(a_n, a_0, \cdots, a_{n-1})~.
\end{equation}
A {\it cyclic cocycle} is a  cyclic cochain $\varphi$ for which $b\varphi=0$.

A straightforward computation shows that the sets of cyclic $n$-cochains  
$\cochala{n}=\{\varphi\in \cocha{n} ~,~ \lambda \varphi = \varphi \} $ are preserved by the Hochschild
boundary operator: $(1-\lambda)\varphi = 0$ implies that $(1-\lambda)b\varphi = 0$. Thus we get a
subcomplex $(C^*_\lambda(\ca)=\bigoplus_n\cochala{n}, b)$ of the complex 
$(C^*(\ca)=\bigoplus_n\cocha{n}, b)$. The  
 {\it cyclic cohomology} $HC^{*}(\ca)$ of the algebra $\ca$ is
the cohomology of this subcomplex,
\begin{equation}
HC^n(\ca) :=  H^n(C^*_\lambda(\ca),b) = Z^n_\lambda / B^n_\lambda~,
\end{equation}
with the cyclic cocycles given by $Z^n_\lambda:=\ker(b :\cochala{n}\ra\cochala{n+1})$ and the
cyclic coboundaries given by
$B^n :=\im(b :\cochala{n-1}\ra\cochala{n})$.

One can also define an operator $B$ which is dual to the one in \eqn{Bi}
for the homology and give a bicomplex description of cyclic cohomology by giving a diagram dual to the
one in \eqn{toco} with all arrows inverted and all indices `up'. Since we shall not need this
description later on, we only refer to \cite{co94} for all details. We mention an additional important
operator, the {\it periodicity operator} $S$ which is a map of degree $2$ between cyclic cocycle,
\begin{align}\label{peop}
& S : Z^{n-1}_\lambda \lra Z^{n+1}_\lambda ~,  \\
& S\varphi(a_0, a_1, \cdots, a_{n+1}) := -\frac{1}{n(n+1)}
\sum_{j=1}^{n} \varphi(a_0, \cdots, a_{j-1}a_j a_{j+1}, \cdots, a_{n+1})  \nn
&  ~~~~~ -\frac{1}{n(n+1)}
\sum_{1\leq i< j\leq n}^{n} (-1)^{i+j} \varphi(a_0, \cdots, a_{i-1}a_i, \cdots, a_ja_{j+1}, \cdots,
a_{n+1})~. \nonumber
\end{align} 
One shows that $S(Z^{n-1}_\lambda) \subseteq Z^{n+1}_\lambda$. In fact $S(Z^{n-1}_\lambda)
\subseteq B^{n+1}$, the latter being the Hochschild coboundary; and cyclicity is easy to
show.\\ The induced morphisms in cohomology
$S:HC^{n} \ra HC^{n+2}$ define two directed systems of abelian groups. Their inductive limits
\begin{equation}\label{pecy}
HP^0(\ca) := \lim_{\ra} HC^{2n}(\ca)~, ~~~HP^1(\ca) := \lim_{\ra} HC^{2n+1}(\ca)~,
\end{equation}
form a $\IZ_2$-graded group which is called the {\it periodic cyclic cohomology} $HP^{*}(\ca)$ of the
algebra $\ca$. \\
`Il va sans dire': there is also a {\it periodic cyclic homology} \cite{co94,lo}.

\subsection{Noncommutative algebras from idempotents}\label{se:ideven}

Let $\ca$ be an algebra (over $\IC$) and let 
$e \in \Mat_r (\ca) \, , \ e^2 = e$, be an idempotent. Its {\it even} (reduced) Chern character is  a
formal  sum of chains 
\begin{equation}
\ChC_*(e) = \sum_k \ChC_k(e)  ~, \label{evenchch}
\end{equation}
with the component $\ChC_k (e)$ an element of
$\ca \ot \ov{\ca}^{\ot 2k} $, 
where $\ov{\ca} = \ca / \IC 1$ is the quotient of $\ca$ by the scalar
multiples of
the unit 1.
The formula for $\ChC_k (e)$ is (with $\lambda_k$ a normalization
constant),
\begin{eqnarray}
\ChC_k (e) &=& \left\lgl \left( e  - \frac{1}{2} \, \II_r
\right) \ot e \ot e \cdots \ot e  \right\rgl \nonumber \\
 &=& \lambda_k \sum \left( e_{i_0 i_1} - \frac{1}{2} \, \delta_{i_0 i_1}
\right) \ot \wt{e}_{i_1 i_2} \ot \wt{e}_{i_2 i_3} \cdots \ot \wt{e}_{i_{2k} i_0} \label{chke}
\end{eqnarray}
where $\delta_{ij}$ is the usual Kronecker symbol and only the class $\wt{e}_{i_j
i_{j+1}} \in \ov{\ca}$ is used in the formula.
The crucial property of the character $\ChC_*(e)$ is that it 
defines a  cycle \cite{co85,co94,lo,co00} in the
reduced $(b,B)$-bicomplex of cyclic homology described above,
\begin{equation} 
(b+B)\ChC_*(e) = 0~, ~~~~ 
B \, \ChC_k (e) = b \, \ChC_{k+1} (e) \, . \label{Bbpro}
\end{equation}
It turns out that the map $e \mapsto \ChC_*(e)$ leads to a well defined map from the $K$ theory group
$K_0(\ca)$  to cyclic homology of $\ca$ (in fact the correct receptacle is period cyclic
homology \cite{lo}). In Section~\ref{se:fmqes}
below, we shall construct some interesting examples of this Chern character on quantum
spheres. For the remaining part of this Section we shall use it to define some `even'
dimensional noncommutative algebras (including spheres).

\bigskip
For any pair of integers $m,r$ we shall construct a universal algebra $\ca_{m,r}$
as follows. We let $\ca_{m,r}$ be generated by the $r^2$ elements $e_{ij}$, 
$i,j \in \{ 1, \ldots , r \}$, $e = [e_{ij}]$ on which we first impose the relations
stating that $e$ is an idempotent 
\begin{equation}
e^2 = e . \label{eq2.6}
\end{equation}
We impose additional relations by requiring the vanishing of all `lower degree' components of
the Chern character of $e$,
\begin{equation}
\ChC_k (e) = 0 ~, \quad \fl \, k < m \, . \label{eq2.5}
\end{equation}
Then, an admissible morphism from $\ca_{m,r}$ to an arbitrary algebra $\cb$,
\begin{equation}
\rho : \ca_{m,r} \ra \cb, 
\end{equation}
is given by the $\rho (e_{ij}) \in \cb$ which fulfill
$\rho (e)^2 = \rho (e)$, and 
\begin{equation}
\ChC_k (\rho(e)) = 0 ~, \quad \fl \, k < m . \label{eq2.7}
\end{equation} We define the algebra 
$\ca_{m,r}$ as the quotient of the algebra defined by (\ref{eq2.6}) by the intersection of kernels of
the admissible morphisms $\rho$. Elements of the algebra $\ca_{m,r}$ can be represented as polynomials
in the generators
$e_{ij}$ and to prove that such a polynomial $P (e_{ij})$ is non zero in
$\ca_{m,r}$ one must construct a solution to the above equations for which $P (e_{ij}) \ne 0$.\\ To get
a \Cs-algebra we endow $\ca_{m,r}$ with the involution given by,
\begin{equation} (e_{ij})^* = e_{ji} \label{eq2.10}
\end{equation} which means that $e = e^*$ in $\Mat_r (\ca)$, i.e. $e$ is a projection in $\Mat_r (\ca)$ 
(or equivalently, a self-adjoint idempotent). We define a norm by,
\begin{equation}
\Vert P \Vert = \Sup \, \Vert (\pi (P)) \Vert \label{eq2.11}
\end{equation} where $\pi$ ranges through all representations of the above equations on Hilbert spaces.
Such a $\pi$ is given by a Hilbert space $\ch$ and a self-adjoint idempotent,
\begin{equation} E \in \Mat_r (\cl (\ch)) \, , \ E^2 = E \, , \ E = E^* \label{eq2.12}
\end{equation} such that (\ref{eq2.7}) holds for $\cb = \cl (\ch)$. For any polynomial $P (e_{ij})$ the
quantity (\ref{eq2.11}), i.e. the supremum of the norms,
$\Vert P (E_{ij}) \Vert $ is finite. \\ We let $A_{m,r}$ be the universal \Cs-algebra obtained as the
completion of $\ca_{m,r}$ for the above norm.

\subsection{Noncommutative algebras from unitaries}\label{se:unodd}

In the odd case, more than projections one rather needs unitary elements and  
the formul{\ae} for the {\it odd } (reduced) Chern character in cyclic homology are
similar to those  above. The Chern  character of a unitary  $u\in\Mat_r(\ca)$ is 
a formal sum of chains
\begin{equation}
\ChC_*(u) = \sum_k \ChC_k(u) ~, \label{oddchch}
\end{equation}
with the component $\ChC_{n+\frac{1}{2}} (u)$ as
element of 
$\ca \ot \ov{\ca}^{\ot(2n-1)}$ given by 
\begin{eqnarray}
 \ChC_{k+\half}(u) &=& \lambda_k \left(u^{i_0}_{i_1} \ot (u^*)^{i_1}_{i_2} 
\ot u^{i_2}_{i_3} \ot \cdots \ot (u^*)^{i_{2k+1}}_{i_0} \right. \nn
&& 
~~~~~~~~~~~~ \left.- (u^*)^{i_0}_{i_1} \ot u^{i_1}_{i_2} 
\ot (u^*)^{i_2}_{i_3} \ot \cdots \ot u^{i_{2k+1}}_{i_0} \right)~, 
\end{eqnarray}
and $\lambda_k$ suitable normalization constants.
Again $\ChC_*(u)$  
defines a  cycle in the
reduced $(b,B)$-bicomplex of cyclic homology \cite{co85,co94,lo,co00},
\begin{equation}
(b+B)\ChC_*(u) = 0~, ~~~~ 
B \, \ChC_{k+\half} (e) = b \, \ChC_{k+\half+1} (e) \, , \label{Bbuni}
\end{equation}
and the map $u \mapsto \ChC_*(u)$ leads to a well defined map from the $K$ theory group $K_1(\ca)$ to
(in fact periodic) cyclic homology.

\bigskip
For any pair of integers $m,r$ we can define $\cb_{m,r}$ to be 
the universal algebra
generated by the $r^2$ elements $u_{ij}$,
$i,j \in \{ 1, \ldots , r \}$, $u = [u_{ij}]$ and we impose as above the relations
\begin{equation}
\ChC_{k+\frac{1}{2}} (\rho(u)) = 0 \qquad \fl \, k < m \, . \label{eq2.7bis}
\end{equation}
To get a  \Cs-algebra we endow $\cb_{m,r}$ with the involution
given by,
\begin{equation}
u \, u^*=u^* \, u = 1 ,\label{eq12.6}
\end{equation}
which means that $u$ is a unitary in $\Mat_r (\ca)$. As before,  we define a norm by,
\begin{equation}
\Vert P \Vert = \Sup \, \Vert (\pi (P)) \Vert 
\end{equation}
where $\pi$ ranges through all representations of the above equations in
Hilbert
space. \\
We let $B_{m,r}$ be the universal \Cs-algebra obtained as the completion
of $\cb_{m,r}$
for the above norm.

\section{Fredholm modules and spectral triples}

As we have mentioned in Section~\ref{inal}, the Chern characters $\ChC_*(x)$ leads to well defined maps
from the $K$ theory groups $K_*(\ca)$  to (period) cyclic homology. The dual 
Chern characters, $\ChC^*$, of even and odd Fredholm modules provides similar maps to (period) cyclic
cohomology.

\subsection{Fredholm modules and index theorems} A Fredholm module can be thought of as an
abstract elliptic operator. The full fledged theory  started with Atiyah and culminated in the
$KK$-theory of Kasparov and the cyclic cohomology of Connes. Here we shall only mention the few facts
that we shall need later on. 

Let  $\ca$ be an algebra  with involution. 
An {\rm odd} Fredholm module \cite{co85} over $\ca$ consists of

\noindent
1) a representation $\psi$ of the algebra $\ca$ on an Hilbert space $\ch$;

\noindent
2) an operator $F$ on $\ch$ such that 
\begin{eqnarray}
&& F^2 = \II ~, ~~~ F^*=F ~,  \nn 
&& [F, \psi(a)] \in \ck  ~, ~~~ \forall ~a\in \ca ~, \label{frmo} 
\end{eqnarray}
where $\ck $ are the compact operators on $\ch$.

\noindent
An {\rm even} Fredholm module has also a $\IZ_2$-grading $\gamma$ of $\ch$, $\gamma^* = \gamma$,
$\gamma^2=\II$, such that  
\begin{eqnarray}
&& F \gamma + \gamma F  = 0 ~, \nn
&& \psi(a) \gamma - \gamma \psi(a) = 0~, ~~~ \forall ~a\in \ca ~. 
\end{eqnarray}
In fact, often the first of conditions \eqn{frmo} needs to be weakened somehow to $F^2 - \II \in
\ck $.\\ With an even module we shall indicate with $\ch^\pm$ and $\psi^\pm$ the component of
the Hilbert space and of the representation with respect to the grading.

Given any positive integer $r$, one can extend the previous modules to a Fredholm module  $(\ch_r,
F_r)$ over the algebra
$\Mat_r(\ca) = \ca \ot \Mat_r(\IC)$ by a simple procedure
\begin{equation}
\ch_r = \ch \ot \IC^r ~, ~~~ \psi_r = \psi \ot \id ~, ~~~ \IF_r = F \ot \II_r ~,
\end{equation}
and $\gamma_r = \gamma \ot \II_r$, for an even module.

\bigskip
The importance of Fredholm modules is testified by the following theorem which can be associated
with the names of Atiyah and Kasparov \cite{at70,ka},
\begin{theorem}~\\
{\rm a)} Let $(\ch, F, \gamma)$ be an even Fredholm module over the algebra $\ca$. And let
$e\in\Mat_r(\ca)$  be a projection $e^2=e=e^*$. Then we have a Fredholm operator
\begin{equation}
\psi^-_r(e) F_r \psi^+_r(e) : \psi^+_r(e) \ch_r \ra \psi^-_r(e) \ch_r ~,
\end{equation}
whose index depends only on the class of the projection $e$ in the $K$-theory of $\ca$. Thus we get 
an additive map
\begin{eqnarray}
&&\varphi : K_0(\ca) \ra \IZ ~, \nn
&&\varphi([e]) = {\rm Index} \left(\psi^-_r(e) F_r \psi^+_r(e)\right)~. \label{evenind}
\end{eqnarray}
{\rm b)} Let $(\ch, F)$ be an odd Fredholm module over the algebra $\ca$, and take the projection 
$E=\frac{1}{2} (\II + F)$. Let $u\in\Mat_r(\ca)$ be 
unitary $u u^* = u^* u = \II$. Then we have a Fredholm operator
\begin{equation}
E_r \psi_r(u) E_r : E_r \ch_r \ra E_r \ch_r ~,
\end{equation} 
whose index depends only on the class of the unitary $u$ in the $K$-theory of $\ca$. Thus we get an
additive map
\begin{eqnarray}
&&\varphi : K_1(\ca) \ra \IZ ~, \nn
&&\varphi([u]) = {\rm Index} \left(E_r \psi_r(u) E_r\right)~. \label{oddind}
\end{eqnarray}
\end{theorem}
If $\ca$ is a \Cs-algebra, both in the even and the odd cases,  the index map $\varphi$
only depends on the {\it $K$-homology} class 
\begin{equation}
[(\ch, F)] \in KK(\ca, \IC) ~,
\end{equation}
of the Fredholm module in the Kasparov $KK$ group, $K^*(\ca)=KK(\ca, \IC)$, which is the abelian 
group of stable homotopy classes of Fredholm modules over $\ca$ \cite{ka}. \\
Both in the even and odd cases, the index pairings \eqn{evenind} and \eqn{oddind} can be given  as
\cite{co94}
\begin{equation}
\varphi(x) = \hs{\ChC^*(\ch, F)}{\ChC_*(x)} ~,  ~~x\in K_*(\ca) ~,
\end{equation}
via the Chern characters
\begin{equation}
\ChC^*(\ch, F) \in HC^*(\ca), ~~~\ChC_*(x) \in HC_*(\ca)~,
\end{equation}
and the pairing between cyclic cohomology $HC^*(\ca)$ and cyclic homology $HC_*(\ca)$ of the algebra
$\ca$. 

The Chern character $\ChC_*(x)$ in homology is given by \eqn{evenchch} and \eqn{oddchch} in the 
even and odd case respectively. As for the Chern character $\ChC^*(x)$ in cohomology we shall  
give some fundamentals in the next Section. 

\subsection{The Chern characters of Fredholm modules}
For the general theory  we refer to \cite{co94}. In
Section~\ref{se:fmqes} we shall construct some interesting examples of these Chern characters on quantum
spheres. Additional examples have been constructed in \cite{tom}.

We recall \cite{si} that on a Hilbert space $\ch$ and with $\ck$ denoting the compact
operators one defines, for
$p\in[1,\infty[$, the Schatten
$p$-class, $\cl^p$, as the ideal of compact operators for which $\Tr T$ is finite: 
$\cl^p =\{T\in\ck ~:~\Tr T < \infty \}$. Then, the H\"older inequality states that 
$\cl^{p_1}\cl^{p_2}\cdots\cl^{p_k}\subset\cl^{p}$, with $p^{-1}=\sum_{j=1}^kp^{-1}_j$. 
 
\bigskip
Let now $(\ch, F)$ be Fredholm module (even or odd) over the algebra $\ca$.
 We say that $(\ch, F)$ is {\it $p$-summable} if
\begin{equation}
[F, \psi(a)] \in \cl^p  ~, ~~~ \forall ~a\in \ca ~.
\end{equation} 
For simplicity, in the rest of this section, we shall drop the symbol $\psi$ which indicates the
representation on $\ca$ on $\ch$. The idea is then to construct `quantized differential forms'
and integrate (via a trace) forms of degree higher enough so that they belong to $\cl^1$. In
fact, one need to introduce a conditional trace. Given an operator $T$ on $\ch$ such that
$FT+TF\in\cl^1$, one defines
\begin{equation}
\Tr' T := \half \Tr F(FT+TF)~; 
\end{equation}
note that, if $T\in\cl^1$ then $\Tr T = \Tr' T$ by cyclicity of the trace.

Let now $n$ be a nonnegative integer and let $(\ch, F)$ be Fredholm module over
the algebra $\ca$. We take this module to be {\it even} or {\it odd} according to whether $n$
is even or odd; and we shall also take it to be $(n+1)$-summable. 
We shall construct a so called $n$-dimensional {\it cycle} $(\Omega^*=\op_k\Omega^k, d, \int
)$ over the algebra $\ca$. Elements of $\Omega^k$ are {\it quantized differential forms}:
$\Omega^0 =\ca$ and for $k>0$, $\Omega^k$ is the linear span of operators of the form
\begin{equation}
\omega = a_0 [F,a_1] \cdots [F,a_n]~, ~~~a_j\in\ca ~.  
\end{equation}
By the assumption of summability, H\"older inequality gives that $\Omega^k\subset
\cl^{\frac{n+1}{k}}$. The product in $\Omega^*$ is just the product of operators $\omega
\omega' \in \Omega^{k+k'}$ for any $\omega\in \Omega^{k}$ and $\omega'\in \Omega^{k'}$.
The differential $d:\Omega^k\ra\Omega^{k+1}$ is defined by
\begin{equation}
d \omega =  F \omega -(-1)^k \omega F ~, ~~~\omega\in \Omega^{k} ~,  
\end{equation}
and $F^2=1$ implies both $d^2=0$ and the fact that $d$ is a graded derivation 
\begin{equation}
d (\omega \omega') = (d \omega) \omega' + (-1)^k \omega d \omega' ~, ~~~  
\omega\in \Omega^{k}~, ~~\omega'\in \Omega^{k'}~. 
\end{equation}
Finally, one defines a trace in degree $n$ by, 
\begin{equation}
\int : \Omega^{n} \ra \IC~,
\end{equation}
which is both closed ($\int d \omega = 0$) and graded ($ \int \omega \omega' 
=  (-1)^{kk'} \int \omega' \omega $). \\
Let us first consider the case $n$ is odd. With $\omega \in \Omega^{n}$ one defines
\begin{equation}\label{intodd}
\int \omega := \Tr' \omega = \half \Tr F(F\omega+\omega F)) = \half \Tr F d \omega ~, 
\end{equation}
which is well defined since $F d \omega\in\cl^1$. \\
If $n$ is even and $\gamma$ is the grading, 
with $\omega \in \Omega^{n}$ one defines
\begin{equation}\label{inteven}
\int \omega := \Tr' \gamma\omega = \half \Tr F(F\gamma\omega+\gamma\omega F)) = 
\half \Tr \gamma F d \omega ~, 
\end{equation}
(remember that $F\gamma=-\gamma F)$; this is again well defined since $\gamma F d
\omega\in\cl^1$.
One straightforwardly proves closeness and graded cyclicity of both the integrals
\eqn{intodd} and \eqn{inteven}.

The {\it character} of the Fredholm module is the cyclic
cocycle $\tau^n \in Z^n_\lambda(\ca)$ given by,
\begin{equation}\label{chfr}
\tau^n(a_0,a_1,\cdots,a_n) := \int a_0 da_1 \cdots da_n ~, ~~~ a_j\in\ca ;
\end{equation}
explicitly, \begin{eqnarray}
&& \tau^n(a_0,a_1,\cdots,a_n) = \Tr' a_0[F,a_1],\cdots,[F,a_n] ~, ~~~n ~~{\rm odd}~, 
\label{chfrodd}\\ 
&& \tau^n(a_0,a_1,\cdots,a_n) = \Tr' \gamma ~a_0[F,a_1],\cdots,[F,a_n] ~, ~~~n ~~{\rm even}~.
\label{chfreven}
\end{eqnarray}
In both cases one checks closure, $b \tau^n = 0$,  and cyclicity, $\lambda \tau^n
= (-1)^n\tau^n$.

We see that there is ambiguity in the choice of the integer $n$. Given a  Fredholm
module $(\ch,F)$ over $\ca$, the parity of $n$ is fixed by for its precise value there is
only a lower bound determined by the $(n+1)$-summability. Indeed, since
$\cl^{p_1}\subset\cl^{p_2}$ if $p_1 \leq p_2$, one could replace $n$ by $n+2k$ with $k$ any
integer. Thus one gets a sequence of cyclic cocycle $\tau^{n+2k} \in Z^{n+2k}_\lambda(\ca),
k\geq 0$, with the same parity. The crucial fact is that the cyclic cohomology classes of
these cocycles are related by the periodicity operator $S$ in \eqn{peop}. The characters
$\tau^{n+2k}$ satisfy 
\begin{equation} 
S [\tau^{m}]_\lambda = c_m [\tau^{m+2}]_\lambda ~, ~~~ {\rm in} ~~~ HC^{m+2}(\ca)~,
~~~m=n+2k~, ~k\geq 0~,
\end{equation}
with $c_m$ a constant depending on $m$ (one could get rid of these constants by suitably
normalizing the characters in \eqn{chfrodd} and \eqn{chfreven}). Therefore, the
sequence $\{\tau^{n+2k} \}_{k\geq 0}$ determine a well defined class $[\tau^{F}]$ in the
periodic cyclic cohomology $HP^0(\ca)$ or $HP^1(\ca)$ according to whether $n$ is even or odd. The
class $[\tau^{F}]$ is the Chern character of the Fredholm module $(\ca, \ch, F)$ in periodic
cyclic cohomology.

\subsection{Spectral triples and index theorems} 

As already mentioned, a noncommutative geometry is described
by a spectral triple
\cite{co94}
\begin{equation}
(\ca , \ch , D) ~.\label{eq1}
\end{equation}
Here $\ca$ is an algebra with involution, together with a representation $\psi$ 
of $\ca$ as bounded operators on a  
Hilbert
space $\ch$ as bounded operators, and $D$ is a self-adjoint operator with compact resolvent
and such that,
\begin{equation}
[D,\psi(a)] \ \hbox{is bounded} \ \forall \, a \in \ca \, . \label{eq2}
\end{equation}
An {\rm even} spectral triple has also a $\IZ_2$-grading $\gamma$ of $\ch$, $\gamma^* = \gamma$,
$\gamma^2=\II$, with the additional properties, 
\begin{eqnarray}
&& D \gamma + \gamma D  = 0 ~, \nn
&& \psi(a)  \gamma - \gamma \psi(a) = 0~, ~~~ \forall ~a\in \ca ~. 
\end{eqnarray}
Given a spectral triple there is associated a fredholm module with the operator $F$ just given by 
the sign of
$D$, $F = D \abs{D}^{-1}$ (if the kernel of $D$ in not trivial one can still adjust things and define
such an $F$).
\\ The operator $D$ plays in general the role of the Dirac operator \cite{lami} in
ordinary Riemannian geometry. It specifies both the $K$-homology fundamental class (cf. \cite{co94}),
as well as the metric on the state space of
$\ca$ by
\begin{equation}
d (\varphi , \psi) = \Sup \, \{ \vert \varphi (a) - \psi (a) \vert ; \Vert [D,a] \Vert
\leq 1 \} ~. \label{eq3}
\end{equation}
What
holds things together
in this spectral point of view on noncommutative geometry is the nontriviality of the pairing
between
the $K$-theory of the algebra $\ca$ and the $K$-homology class of $D$. There are index maps as  with
Fredholm modules above,  
\begin{equation}
\varphi : K_* (\ca) \ra \IZ \,  \label{eq4}
\end{equation}
and the maps $\varphi$ given by expressions like \eqn{evenind} and \eqn{oddind} with the operator 
$D$ replacing the operator $F$ there. 
 
An operator theoretic index formula
\cite{co94}, \cite{como}, \cite{grvafi} expresses the above index
pairing (\ref{eq4}) by explicit
{\it local} cyclic cocycles on the algebra $\ca$. These local formulas become
extremely simple in the special case where only the top component of the Chern
character $\ChC_*(e)$ in cyclic homology fails to vanish. This is easy
to understand in the
analogous simpler case
of ordinary manifolds since the Atiyah-Singer index formula gives the
integral of
the product of the Chern character $\ChC (E)$, of the bundle $E$ over
the manifold
$M$, by the index class; if the only component of
$\ChC (E)$ is $\ChC_n$, $n = \frac{1}{2} \, \dim M$ only the
$0$-dimensional
component of the index class is involved in the index formula.

For instance, in the even case, provided the components $\ChC_k (e)$ all vanish for $k < n$ the 
index formula reduces to the following,
\begin{equation}
\varphi(e) = (-1)^n \int\!\!\!\!\!\!- \ \gamma \left( e - \frac{1}{2}
\right) [D,e]^{2n} \, D^{-2n} ~.\label{eq8}
\end{equation}
Here, $e$ is a projection $e^2=e=e^*$,  $\gamma$
is the
$\IZ / 2$ grading of $\ch$ as above, the resolvent of $D$ is of order
$\frac{1}{2n}$ (i.e. its characteristic values $\mu_k$ are
$0(k^{-\frac{1}{2n}})$)
and $\int\!\!\!\!\!-$ is the coefficient of the logarithmic divergency in the
ordinary operator trace \cite{di} \cite{wod}. There is a similar formula for the odd case.

\begin{example}\label{ex:ctm}{The Canonical Triple over a Manifold}\\
The basic example of spectral triple is the {\it canonical triple} on a
closed $n$-dimensional Riemannian spin manifold $(M,g)$. A spin manifold is a manifold on which it is
possible to construct principal bundles having the groups $Spin(n)$ as structure groups. A manifold
admits a spin structure if and only if its second Stiefel-Whitney class vanishes \cite{lami}. 

The canonical spectral triple $(\ca, \ch, D)$ over the manifold $M$ is as follows:\\
$1)$~~$\ca=\cinf(M)$ is the algebra of complex valued smooth functions on $M$. \\
$2)$~~$\ch=L^2(M,\cs)$ is the Hilbert space of square integrable sections of
the irreducible, rank $2^{[n/2]}$, spinor bundle over $M$ ; its elements
are spinor fields over $M$.   The scalar product in
$L^2(M,\cs)$ is the usual one of the measure
$d\mu(g)$ of the metric $g$, $(\psi, \phi) = \int d\mu(g) \bar{\psi(x)}\cdot \phi(x)$,   
with the pointwise scalar product in the spinor space being
the natural one in $\IC^{2^{[n/2]}}$. \\
$3)$~~$D$ is the Dirac operator of the Levi-Civita
connection of the metric $g$. 
It can be written locally as 
\begin{equation}
D = \gamma^\mu(x)(\pa_\mu+\omega_\mu^\cs) ~, 
\label{dirac}
\end{equation} 
where $\omega_\mu^\cs$ is 
the lift of the Levi-Civita connection to the bundle of spinors. 
The curved gamma matrices $\{\gamma^\mu(x)\}$ are 
Hermitian and satisfy
\begin{equation}\label{gamrel} 
\gamma^\mu(x)\gamma^\nu(x) + \gamma^\nu(x)\gamma^\mu(x) = 2g(dx^\mu, dx^n) =
2g^{\mu\nu}~, ~~\mu, \nu = 1, \dots, n~.
\end{equation}
The elements of the algebra $\ca$ act as multiplicative
operators on
$\ch$, 
\begin{equation}(f \psi)(x)=:f(x)\psi(x)~, ~~\forall~ f\in\ca~, \psi\in\ch~. \label{mul}
\end{equation}
\end{example}
For this triple, the distance in \eqn{eq3} is the geodesic distance on the manifold $M$ of 
the metric $g$.

An additional important ingredient is provided by a {\it real structure}. In the
context of the canonical triple, this is given by $J$, the charge conjugation
operator, which is an antilinear isometry of $\ch$. We refer to \cite{co94} for all details; 
a friendly introduction is in \cite{la97}.

\section{Examples of Isospectral Deformations}

We shall now construct some examples of (a priori noncommutative) spaces ${\rm Gr}_{m,r}$ such
that
\begin{equation}
A_{m,r} = C ({\rm Gr}_{m,r}) \qquad {\rm or } \qquad B_{m,r} = C ({\rm Gr}_{m,r}) ,\label{eq2.14}
\end{equation}
according to even or odd dimensions, with the $\Cs$-algebras
$A_{m,r}$ and $B_{m,r}$ defined at the end of Section~\ref{se:ideven} and Section~\ref{se:unodd}, and
associated with the vanishing of the `lower degree' components of the Chern character of an idempotent
and of a unitary respectively.

\subsection{Spheres in dimension $2$}
The simplest case is $m=1$, $r=2$. We have then
\begin{equation}e =  
\begin{pmatrix} e_{11} &e_{12} \cr e_{21} &e_{22}
\end{pmatrix} 
\end{equation}
and the condition (\ref{eq2.7}) just means that
\begin{equation}
e_{11} + e_{22} = 1 \label{eq2.15}
\end{equation}
while (\ref{eq2.6}) means that
\begin{eqnarray}
&e_{11}^2 + e_{12} \, e_{21} = e_{11} \, , \ e_{11} \, e_{12} + e_{12} \,
e_{22} =
e_{12} \, , \label{eq2.16} \\
&e_{21} \, e_{11} + e_{22} \, e_{21} = e_{21} \, , \ e_{21} \, e_{12} +
e_{22}^2 =
e_{22} \, . \nonumber
\end{eqnarray}
By (\ref{eq2.15}) we get $e_{11} - e_{11}^2 = e_{22} - e_{22}^2$,
so that
(\ref{eq2.16}) shows that $e_{12} \, e_{21} = e_{21} \, e_{12}$. We also
see that
$e_{12}$ and $e_{21}$ both commute with $e_{11}$. This shows that
$\ca_{1,2}$ is
commutative and allows to check that ${\rm Gr}_{1,2} = S^2$ is the
2-sphere. Thus
${\rm Gr}_{1,2}$ is an ordinary commutative space.

\subsection{Spheres in dimension $4$}\label{sphd4}
Next, we move on to the case $m=2$, $r=4$.  

Note first that the notion of admissible 
morphism is a non trivial piece of structure on ${\rm Gr}_{2,4}$
since, for instance, the identity map is not admissible \cite{codu-vi}.

Commutative solutions were found in \cite{co00} with the commutative algebra $A=C(S^4)$ and an admissible
surjection $A_{2,4} \ra C(S^4)$, where the sphere $S^4$ appears naturally as quaternionic projective
space, $S^4 = P_1(\IH)$. 

In \cite{cola} we found noncommutative solutions, showing that the algebra $A_{2,4}$ is
noncommutative, and we constructed explicit admissible surjections,
\begin{equation}
A_{2,4} \ra C(S_{\theta}^4)
\end{equation}
where
$S_{\theta}^4$ is the noncommutative 4-sphere we are about to describe and whose form is dictated by
natural deformations of the ordinary $4$-sphere, similar in spirit to the standard deformation of the
torus $\IT^2$ to the noncommutative torus $\IT^2_\theta$. In fact, as will become evident later on,
noncommutative tori in arbitrary dimensions play a central role in the deformations.
\\ We first determine the algebra generated by the usual matrices $\Mat_4(\IC)$ and a projection $e = e^*
= e^2$ such that
$\ChC_0 (e) = 0$ as above and whose matrix expression is of the form,
\begin{equation}
[ e^{ij} ] = \half  \begin{pmatrix} 
 q_{11} &q_{12} \cr q_{21} &q_{22}  
\end{pmatrix}
\label{eq58}
\end{equation}
where each $q_{ij}$ is a $2 \ts 2$ matrix of the form,
\begin{equation}
q =  \begin{pmatrix}
\alpha &\beta \cr -\lambda \beta^* &\alpha^* 
\end{pmatrix}
\, ,
\label{eq59}
\end{equation}
and $\lambda=\exp (2\pi i \theta)$ is a complex number of modulus one 
(different from $-1$ for convenience).
Since $e = e^*$, both $q_{11}$ and $q_{22}$ are self-adjoint, moreover since
 $\ChC_0 (e) = 0$, we can
find $z = z^*$ such that,
\begin{equation}
q_{11} =  \begin{pmatrix}
1+z &0 \cr 0 &1+z   \end{pmatrix} \, , \
q_{22} = 
\begin{pmatrix} 1-z &0 \cr 0 &1-z  \end{pmatrix} \, . \label{eq60}
\end{equation}
We let $q_{12} =
\left( \begin{smallmatrix}
\alpha  & \beta\\
-\lambda \beta^*  & \alpha^*
\end{smallmatrix}\right) $,
we then get from $e = e^*$,
\begin{equation}
q_{21} =   \begin{pmatrix} \alpha^* &- \bar{\lambda} \beta\cr \beta^* &\alpha  
\end{pmatrix}
\, . \label{eq61}
\end{equation}
We thus see that the commutant $\ca_{\theta}$ of $\Mat_4 (\IC)$ is generated
by $z,\alpha,\beta$
and we first need to find the relations imposed by the equality $e^2 = e$.
In terms of the matrix 
\begin{equation}
e = \half  \begin{pmatrix}  1+z &q \cr q^* &1-z \end{pmatrix} \, ,
\label{idempot}
\end{equation}
the equation $e^2 = e$ means that $z^2 + q q^* = 1$, $z^2 + q^* q =
1$ and $[z,q] = 0$. This shows that $z$ commutes with $\alpha$, $\beta$, $\alpha^*$ and
$\beta^*$ and since $qq^* = q^* q$ is a diagonal matrix
\begin{equation}
\alpha \alpha^* = \alpha^* \alpha \, , \ \alpha \beta=\lambda  \beta\alpha \, , \ \alpha^* \beta
=\bar{\lambda} \beta\alpha^* \, , \ \beta\beta^* = \beta^* \beta\label{eq62}
\end{equation}
so that the generated algebra $\ca_{\theta}$ is not commutative for $
\lambda$
different from
1. The only further relation, besides $z = z^*$, is a sphere relation
\begin{equation}
\alpha \alpha^* + \beta\beta^* + z^2 = 1 \, . \label{eq63}
\end{equation}
We denote by $S^4_{\theta}$ the corresponding noncommutative space defined by `duality',
so that its algebra of polynomial functions is $\ca(S^4_{\theta})=\ca_{\theta}$.  
This algebra is a deformation of the commutative $*$-algebra $\ca(S^4)$ of
complex polynomial functions on the usual sphere $S^4$ to which it reduces for
$\theta=0$. 

The projection $e$ given in \eqn{idempot} is clearly an element in the matrix algebra
$\Mat_4(\ca_{\theta})\simeq\Mat_4(\IC)\ot\ca_{\theta}$. Then, it naturally acts on the free
$\ca_{\theta}$-module $\ca_{\theta}^4\simeq\IC^4\ot\ca_{\theta}$ and one gets as its range a finite
projective module which can be thought of as the module of `section of a vector bundle' over
$S^4_{\theta}$. The module $e \ca_{\theta}^4$ is a deformation of the usual \cite{at79} complex rank
$2$ instanton bundle over $S^4$ to which it reduces for $\theta=0$ \cite{la00}. 

For the sphere $S^4_{\theta}$ the deformed instanton has correct characteristic classes. The fact
that $\ChC_0(e)$ has been imposed from the very beginning and could be interpreted as stating the
fact that the projection and the corresponding module (the `vector bundle') has complex rank equal
to $2$. Next, we shall check that  the two
dimensional component $\ChC_1(e)$ of the Chern character, automatically vanishes as an element of the
(reduced) $(b,B)$-bicomplex. \\
With
$q = \left( \begin{smallmatrix}
\alpha  & \beta\\
-\lambda \beta^*  & \alpha^*
\end{smallmatrix}\right)$,
we get,
$$
\ChC_1(e) =  \frac{1}{2^3} \biggl\lgl z
\, (dq \, dq^* - dq^* \, dq)  
 + \, q \, (dq^* \, dz - dz \, dq^*) + q^* \, (dz \, dq - dq \, dz )
\biggl\rgl  
$$
where the expectation in the right hand side is relative to $\Mat_2 (\IC)$
(it is a partial trace) and
we use
the notation $\, d \, $ instead of the tensor notation. The diagonal elements of  $\omega = dq \,
dq^*$ are
$$
\omega_{11} = d \alpha \, d\alpha^* + d\beta \, d\beta^* \, , \ \omega_{22} = d\beta^* \, d\beta +
d\alpha^* \, d\alpha
$$
while for $\omega' = dq^* \, dq$ we get,
$$
\omega'_{11} = d\alpha^* \, d\alpha + d\beta \, d\beta^* \, , \ \omega'_{22} = d\beta^* \, d\beta +
d\alpha
\, d\alpha^* \, .
$$
It follows that, since $z$ is diagonal,
\begin{equation}
\biggl\lgl z \, (dq \, dq^* - dq^* \, dq)
\biggl\rgl = 0 \, . \label{eq70}
\end{equation}
The diagonal elements of $q \, dq^* \, dz = \rho$ are
$$
\rho_{11} = \alpha \, d\alpha^* \, dz + \beta \, d\beta^* \, dz \, , \ \rho_{22} = \beta^* \,
d\beta \, dz + \alpha^* \, d\alpha \, dz
$$
while for $\rho' = q^* \, dq \, dz$ they are
$$
\rho'_{11} = \alpha^* \, d\alpha \, dz + \beta \, d\beta^* \, dz \, , \ \rho'_{22} = \beta^* \,
d\beta \, dz + \alpha \, d\alpha^* \, dz \, .
$$
Similarly for $\sigma = q \, dz \, dq^*$ and $\sigma' = q^* \, dz \, dq$ one gets the
required cancellations so that,
\begin{equation}
  \ChC_1(e)  = 0 \, ,
\label{eq71}
\end{equation}
Summing up we thus get that the element  $e \in C^{\infty} (S^4_{\theta}, \Mat_4 (\IC))$ given in 
\eqref{idempot} is a self-adjoint idempotent, $e = e^2 = e^*$, 
and satisfies $\ChC_k (e) = 0$ $\fl \, k < 2$.
Moreover, ${\rm Gr}_{2,4}$ is a noncommutative space and
$S^4_{\theta} \sbs {\rm Gr}_{2,4} $. \\
Since $ \ChC_1(e)=0$, it follows that $ \ChC_2(e)$ is a Hochschild
cycle which will play the role of the
round volume form on $S^4_{\theta}$ and that we shall now compute.
With the above notations one has,
\begin{equation}
  \ChC_2(e)
  = \frac{1}{2^5} \left\lgl   \begin{pmatrix}  z &q \cr q^* & -z 
\end{pmatrix}
\begin{pmatrix}  dz  &dq \cr dq^* & -dz  
\end{pmatrix}^4\right\rgl \, .
\label{ch4}
\end{equation}
The direct computation gives the Hochschild cycle $\ChC_2(e) $ as a sum of five components
\begin{equation}
\ChC_2(e) = z \, c_z + \alpha \, c_{\alpha} + \alpha^* \, c_{\alpha^*} + \beta \, c_{\beta} +
\beta^* \, c_{\beta^*} \, ; \label{hocy}
\end{equation}
where the components $c_z \, , c_{\alpha} \, , c_{\alpha^*} \, , c_{\beta} \, ,
c_{\beta^*}$, which are
elements in the four-fold tensor
product $\ca_{\theta}\ot\ov{\ca_{\theta}}\ot\ov{\ca_{\theta}}\ot\ov{\ca_{\theta}}$, are explicitly
given in \cite{cola}. The vanishing of $b \ChC_2(e)$, which has six hundred terms,
can be checked directly from the commutation relations (\ref{eq62}). The
cycle $\ChC_2(e)$ is totally
`$\lambda$-antisymmetric'.

\bigskip

Our sphere $S^4_{\theta}$ is by construction the suspension of the
noncommutative 3-sphere $S^3_{\theta}$ whose coordinate algebra is
generated by $\alpha$ and $\beta$ as above and say the special value $z=0$.
This 3-sphere is part of a family of spheres that we shall describe in
the next Section. 

Had we taken the deformation parameter to be real, $\lambda=q\in
\IR$, the corresponding 3-sphere $S^3_{q}$ would coincide
with the quantum group $SU(2)_q$.
Similarly, had we taken the deformation parameter in $S^4_{\theta}$ to be real like in
\cite{dalama} we would have obtained a different deformation $S^4_{q}$ of the
commutative
sphere $S^4$, whose algebra is different from the above one. More
important, the  component $ \ChC_1(e)$ of the Chern character
would not vanish  \cite{dala}.

\subsection{Spheres in dimension $3$}
Odd dimensional spaces, in particular spheres, are constructed out of
unitaries rather than projections \cite{cola,codu-vi,codu-vi1}.

Let us consider the lowest dimensional case for which $m=2, r=2$. 
We shall use the convention that repeated indices are summed on. Greek
indices like $\mu,\nu, \cdots, $ are taken to be valued in
$\{0,1,2,3\}$ while latin indices like $j, k, \cdots, $ are taken
to be valued in $\{1,2,3\}$.

We are then looking for an algebra $\cb$ such that \\
$1)$~~$\cb$ is generated as a unital
$*$-algebra by the entries of a unitary matrix 
\begin{equation}
u\in \Mat_2(\cb) \simeq
\Mat_2(\IC) \ot \cb, ~~~u u^* = u^* u = 1, 
\end{equation}
$2)$~~the unitary $u$ satisfies the additional condition
\begin{equation}
\ChC_{\frac{1}{2}}(u) := \sum u^j_i \ot (u^*)^j_i - (u^*)^j_i \ot u^j_i =
0~. \label{chhalf}
\end{equation}
Let us take as `generators' of $\cb$ elements $z^\mu, z^{\mu *}$, 
$\mu \in \{0,1,2,3\}$. Then using ordinary Pauli matrices $\sigma_k$, 
$k \in \{1,2,3\}$, an element in $u\in \Mat_2(\cb)$ can be written as 
\begin{equation}
u = \II_2 z^0 + \sigma_k z^k ~. 
\end{equation}
The requirement that $u$ be unitaries give the following conditions on the
generators
\begin{eqnarray}
&& z^k z^{0 *} - z^0 z^{k *} + \varepsilon_{klm}z^{l} z^{m *} = 0 ~, \nn
&& z^{0 *} z^k - z^{k *} z^0 + \varepsilon_{klm}z^{l *} z^{m } = 0 ~,\nn
&& \sum_{\mu=0}^3 (z^{\mu} z^{\mu *} - z^{\mu *} z^{\mu }) = 0~,
\label{com3sp}
\end{eqnarray}
together with the condition that 
\begin{equation}
\sum_{\mu=0}^3 z^{\mu *} z^{\nu } = 1~. \label{sph3d}
\end{equation}
Notice that the `sphere' relation \eqn{sph3d} is consistent with the
relations \eqn{com3sp} since the latter imply that $\sum_{\mu=0}^3
z^{\mu *} z^{\nu }$ is in the center of $\cb$.\\
Then, one imposes condition \eqn{chhalf} which reads
\begin{equation}
\sum_{\mu=0}^3 (z^{\mu *} \ot z^{\mu} - z^{\mu} \ot z^{\mu *}) = 0~,
\end{equation}
and which is satisfied \cite{codu-vi,codu-vi1} if and only if there exists a
symmetric unitary matrix $\Lambda \in \Mat_4(\IC)$ such that 
\begin{equation}
z^{\mu *} = \Lambda^{\mu}_{\nu} z^{\nu}~. \label{lambda}
\end{equation}
Now, there is some freedom in the definition of the algebra $\cb$ which
is stated by the fact that the defining conditions $1.)$ and $2.)$ above
do not change if we transform 
\begin{equation}
z^{\mu} \mapsto \rho S^{\mu}_{\nu} z^{\nu}~,
\end{equation}
with $\rho \in \U(1)$ and $S\in \SO(4)$. 
Under this transformation, the matrix $\Lambda$ in \eqn{lambda}
transforms as 
\begin{equation} 
\Lambda \mapsto \rho^2 S^t \Lambda S ~. 
\end{equation}
One can then diagonalize the symmetric unitary $\Lambda$ by a real
rotation $S$ and fix its first eigenvalue to be $1$ by an appropriate
choice of $\rho\in\U(1)$. So, we can take
\begin{equation}
\Lambda = \diag(1, e^{- i \varphi_1},e^{- i \varphi_2},e^{- i \varphi_3})~,
\end{equation}
that is, we can put
\begin{equation}
z^0 = x^0, ~~z^k = e^{i \varphi_k} x^k, ~~k\in\{1,2,3\}~,
\end{equation}
with $e^{- i \varphi_k}\in\U(1)$ and $(x^\mu)^*=x^\mu$. 
Conditions \eqn{com3sp} translate to
\begin{eqnarray}
&& 
[x^0, x^k]_{-} \cos \varphi_k = \ii ~[x^l, x^m]_{+} \sin(\varphi_l -\varphi_m)~,
\nn
&& 
[x^0, x^k]_{+} \sin \varphi_k = \ii ~[x^l, x^m]_{-} \cos(\varphi_l -\varphi_m)~, 
\end{eqnarray}
with $(k,l,m)$ the cyclic permutation of $(1,2,3)$ starting with $k=1,2,3$ and 
$[x,y]_{\pm}=xy-yx$. There is also the sphere relation \eqn{sph3d},
\begin{equation}
\sum_{\mu=0}^3 (x^{\mu *})^2 = 1~. \label{sph3dbis}
\end{equation}  
We have therefore a three parameters family of algebras 
$\cb_{{\bf \varphi}}$ which are labelled by an element 
${\bf \varphi} = (e^{- i \varphi_1}, e^{- i \varphi_2}, e^{- i \varphi_3}) 
\in \IT^3$. The algebras $\cb_{{\bf \varphi}}$ are deformations of the
algebra $A(S^3)$ of polynomial functions on an ordinary
$3$-sphere $S^3$
which is obtained for the special value ${\bf \varphi}=(1,1,1)$. We denote
by $S^3_{{\bf \varphi}}$ the corresponding noncommutative space, so
that $A(S^3_{{\bf \varphi}})=\cb_{{\bf \varphi}}$. 
Next, one computes $\ChC_{\frac{3}{2}}(u_{{\bf \varphi}}))$ and shows that
is a non trivial cycle ($b \ChC_{\frac{3}{2}}(u_{{\bf \varphi}})) = 0$)
on $\cb_{{\bf \varphi}}$ \cite{codu-vi}. 

A special value of the parameter $\varphi$ gives the
$3$-sphere $S^3_{\theta}$ described at the end of previous Section. Indeed, 
put $\varphi_1=\varphi_2=-\pi \theta$ and $\varphi_3=0$ and define
\begin{eqnarray}
&& \alpha=x^0 + \ii x^3 ~, ~~\alpha^*=x^0 - \ii x^3 ~, \nn
&& \beta=x^1 + \ii x^2 ~, ~~\beta^*=x^1 - \ii x^2~.
\end{eqnarray}
then $\alpha, \alpha^*, \beta, \beta^*$ satisfies conditions \eqn{eq62},
with $\lambda = \exp(2 \pi \ii \theta)$, together with the relation $\alpha \alpha^* + \beta
\beta^* = 1$, thus defining the sphere $S^3_\theta$ of
Section~\ref{sphd4}.

In Section~\ref{ncspman} we shall describe some higher dimensional
examples.

\subsection{The noncommutative geometry of $S^4_{\theta}$}

Next we will analyze the metric structure, via a Dirac operator
$D$, on our noncommutative 4-spheres $S_{\theta}^4$. 
The operator $D$ will give a solution to the following quartic equation,
\begin{equation}
\left\langle \left( e - \frac{1}{2} \right) [D,e]^4 \right\rangle = \gamma
\label{eq11}
\end{equation}
where $\langle \ \rangle$ is the projection on the commutant of $4 \ts 4$ 
$\IC$-matrices (in fact, it is a partial trace on the matrix entries) and 
$\gamma=\gamma_5$, in the present four
dimensional case, is the grading operator.

Let $C^{\infty} (S_{\theta}^4)$ be the algebra of smooth functions on the
noncommutative sphere $S_{\theta}^4$. We shall construct a spectral triple
$(C^{\infty} (S_{\theta}^4), \ch , D)$
which describes the geometry on $S_{\theta}^4$ corresponding to the round metric. \\
In order to do that we
first need to find good coordinates on
$S_{\theta}^4$ in terms of
which the operator $D$ will be easily expressed. We choose to parametrize
$\alpha , \beta$ and
$z$ as follows,
\begin{equation}
\alpha = \,  u \, \cos \varphi \cos \psi  \ , \ \beta = \, v \,
\sin \varphi \cos
\psi  \, , \ z =  \sin \psi \, . \label{eq2bis}
\end{equation}
Here $\varphi$ and $\psi$ are ordinary angles with domain
$0 \leq \varphi \leq \frac{\pi}{2} \, , \ - \frac{\pi}{2} \leq \psi \leq
\frac{\pi}{2}$,
while $u$ and $v$ are the usual unitary generators of the algebra $C^{\infty}
(\IT_{\theta}^2)$ of smooth functions on the noncommutative 2-torus. Thus
the presentation
of their relations is
\begin{equation}
uv = \lambda v u \, , \ uu^* = u^* u = 1 \, , \ vv^* = v^* v = 1 \, .
\label{eq4bis}
\end{equation}
One checks that $\alpha , \beta , z$ given by (2) satisfy the basic presentation
of the
generators of $C^{\infty} (S_{\theta}^4)$ which thus appears as a {\it
subalgebra} of the
algebra generated (and then closed under smooth calculus) by $e^{i\varphi}$,
$e^{i\psi}$,
$u$ and $v$.

 For $\theta = 0$ the round metric is given as,
\begin{equation}
G = d\alpha \, d \ov \alpha + d \beta \, d \ov \beta + dz^2  \label{eq5bis}
\end{equation}
and in terms of the coordinates, $\varphi , \psi , u , v$ one gets,
\begin{equation}
G = \cos^2 \varphi \cos^2 \psi \, du \, d \ov u + \sin^2 \varphi \cos^2 \psi \, dv
\, d \ov v +
\cos^2 \psi \, d \varphi^2 + d \psi^2 \, . \label{eq6bis}
\end{equation}
Its volume form is given by
\begin{equation}
\omega = \frac{1}{2} \sin \varphi \cos \varphi \, (\cos \psi)^3 \, \ov u \, du \wedge \ov v \, dv
\wedge d\varphi \wedge d\psi \, . \label{eq7bis}
\end{equation}
In terms of these rectangular coordinates we get the following simple
expression for
the Dirac operator,
\begin{eqnarray}
D &= &(\cos \varphi \cos \psi)^{-1} \, u \, \frac{\partial}{\partial u} \, \gamma_1
+ (\sin \varphi
\cos \psi)^{-1} \, v \, \frac{\partial}{\partial v} \, \gamma_2 +  \\
&+ &\frac{ \ii }{\cos \psi} \, (\frac{\partial}{\partial \varphi} \, +
\frac{1}{2} \, {\rm
cotg} \, \varphi - \frac{1}{2} \,
{\rm tg} \, \varphi) \,
\gamma_3 + \ii  \, (\frac{\partial}{\partial \psi} \, - \, \frac{3}{2} \, {\rm tg} \,
\psi) \, \gamma_4 \, . \nonumber
\end{eqnarray}
Here $\gamma_{\mu}$ are the usual Dirac $4 \ts 4$ matrices with
\begin{equation}
\{ \gamma_{\mu} , \gamma_{\nu} \} = 2 \, \delta_{\mu \nu} \, , \ \gamma_{\mu}^* = \gamma_{\mu}
\, .
\label{eq9ter}
\end{equation}
 It is now easy to move on to the noncommutative case, the only
tricky point
is that there are nontrivial boundary conditions for the operator $D$, which are 
in particular antiperiodic in the arguments of both $u$ and $v$.
We shall just
leave them unchanged in the noncommutative case, the only thing which changes is the
algebra and the way it acts in
the Hilbert space as we shall explain in more detail in the next section.
The formula for the operator $D$ is now,
\begin{eqnarray}
D &= &(\cos \varphi \cos \psi)^{-1} \,  \delta_1 \, \gamma_1 + (\sin \varphi
\cos \psi)^{-1} \, \delta_2  \, \gamma_2 +  \\
&+ &\frac{\ii}{\cos \psi}  \, (\frac{\partial}{\partial \varphi} \, +
\frac{1}{2} \, {\rm cotg} \, \varphi - \frac{1}{2} \,
{\rm tg} \, \varphi) \,
\gamma_3  +  \ii  \, (\frac{\partial}{\partial \psi} \, - \, \frac{3}{2} \,
{\rm tg} \,
\psi) \, \gamma_4 \, . \nonumber
\end{eqnarray}
where the $\gamma_\mu$ are the usual Dirac matrices and where $\delta_1$ and
$\delta_2$ are the derivations of the noncommutative torus so that
\begin{eqnarray}
&& \delta_1(u) = u \, , \, \, \delta_1(v) = 0 \, ,  \\
&& \delta_2(u) = 0 \, , \, \, \delta_2(v) = v \, ; \nonumber
\end{eqnarray}
One can then check that the corresponding metric is the round one.

 In order to compute the operator $\left\langle \left( e -
\frac{1}{2} \right) [D,e]^{4} \right\rangle$
(in the tensor product by $\Mat_4( \IC)$) we need
the commutators of $D$ with the generators of $C^{\infty} (S^{4}_{\theta})$.
They are
given by the
following simple expressions,
\begin{eqnarray}
&&[ D, \alpha ] = u \, \{\gamma_1 - \ii \sin\phi \, \gamma_3 - \ii \cos\phi 
\sin\psi \, \gamma_4 \, \} \, , \\
&& [ D, \alpha^* ] = - u^* \, \{\gamma_1 + \ii \sin\phi \, \gamma_3 + \ii \cos\phi
\sin\psi \, \gamma_4 \, \} \, ,  \nonumber \\
&& [ D, \beta ] = v \,  \{\gamma_2 + \ii \cos\phi  \, \gamma_3 - \ii \sin\phi
\sin\psi \, \gamma_4 \, \} \, ,
\nonumber \\
&&  [ D, \beta^* ] = - v^* \, \{\gamma_2 - \ii \cos\phi  \, \gamma_3 +  
\ii \sin\phi \sin\psi  \, \gamma_4 \, \} \, , \nonumber \\
&& [ D, z ] = \ii \cos\psi \, \gamma_4 \, . \nonumber
\end{eqnarray} We check in particular that they are all bounded operators and
hence that
for any $f \in C^{\infty} (S^{4}_{\theta})$ the commutator $[D, \, f]$ is
bounded.
Then, a long but straightforward calculation shows that equation \eqn{eq11} is valid: the operator
$\left\langle \left( e - \frac{1}{2} \right) [D,e]^{4} \right\rangle$ is
a multiple of $\gamma =\gamma_5 := \gamma_1 \gamma_2 \gamma_3 \gamma_4$.
One first checks that it is equal to $\pi (\ChC_2(e))$ where $\ChC_2(e)$ is the Hochschild
cycle in \eqn{hocy}
  and $\pi$ is the canonical map from the
Hochschild chains to operators given by
\begin{equation}
\pi( a_0 \ot a_1 \ot ...\ot a_n) = a_0 [D, a_1]...[D, a_n] \, . \label{eq15bis}
\end{equation}

\subsection{Isospectral noncommutative geometries}
We shall describe fully a noncommutative geometry for $S^{4}_{\theta}$ with the couple 
$(\ch, D)$ just the `commutative' ones associated with the commutative sphere $S^{4}$; hence
realizing an isospectral deformation. 
We shall in fact describe a very general construction of isospectral
deformations
of noncommutative geometries
which implies in particular that any
  compact spin Riemannian manifold $M$ whose isometry group has rank
$\geq 2$ admits a
natural one-parameter isospectral deformation to noncommutative geometries
$M_\theta$. The deformation of the algebra will be performed
along the lines of \cite{ri93}  (see also \cite{va} and \cite{sit1}).

Let us start with the canonical spectral triple $(\ca=C^{\infty}(S^{4}), \ch , D)$ associated with 
the sphere $S^{4}$. We recall  that $\ch= L^2(S^{4},\cs)$ is the Hilbert space of spinors and $D$ is the
Dirac  operator. Also, there is a real structure provided by $J$, the charge conjugation operator, which
is an antilinear isometry of $\ch$.

Recall that on the sphere $S^{4}$ there is an isometric action of the $2$-torus, 
$ \IT^2 \subset {\rm Isom}(S^{4})$ with $\IT = \IR / 2 \pi \IZ$ the usual torus. We let $U(s) , s
\in \IT^2$, be
the corresponding (projective) unitary representation in $\ch = L^2(S^{4},\cs)$ so that by construction
\begin{equation}
U(s) \, D = D \, U(s) \, , \, \, \, U(s) \, J = J \, U(s) \, .
\end{equation}
Also,
\begin{equation}
U(s) \, a \, U(s)^{-1} = \alpha_s(a) \, , \, \, \, \fl \, a \in \ca \, ,
\label{actfun}
\end{equation}
where $\alpha_s \in {\rm Aut}(\ca)$ is the action by isometries on functions on
$S^{4}$. \\
We let $p = (p_1, p_2)$ be the generator of the two-parameter group $U(s)$
so that
\begin{equation}
U(s) = \exp(2 \pi i(s_1 p_1 + s_2 p_2)) \, .
\end{equation}
The operators $p_1$ and $p_2$ commute with $D$ but anticommute with $J$ (due to the antilinearity 
of the latter). Both $p_1$ and $p_2$
have half-integral spectrum,
\begin{equation}
{\rm Spec}(2 \, p_j) \subset \IZ \, , \, \, j = 1, 2 \, .
\end{equation}
Next, we define a bigrading of the algebra of bounded operators in $\ch$ with the
operator $T$ declared to be of bidegree
$(n_1,n_2)$ when,
\begin{equation}
\alpha_s(T) := U(s) \, T \, U(s)^{-1} = \exp(2 \pi i(s_1 n_1 + s_2 n_2)) \, T \, , \, \, \, \fl \, 
s \in
\IT^2 \, . \label{t2act}
\end{equation}
Any operator $T$ of class $C^\infty$ relative to $\alpha_s$ (i. e. such that
the map $s \rightarrow \alpha_s(T) $ is of class $C^\infty$ for the
norm topology) can be uniquely
written as a doubly infinite
norm convergent sum of homogeneous elements,
\begin{equation}
T = \sum_{n_1,n_2} \, \widehat{T}_{n_1,n_2} \, ,
\end{equation}
with $\widehat{T}_{n_1,n_2}$ of bidegree $(n_1,n_2)$ and where the sequence
of norms $||
\widehat{T}_{n_1,n_2} ||$ is of
rapid decay in $(n_1,n_2)$.

Let now $\lambda = \exp(2 \pi \ii \theta)$. For any operator $T$ in $\ch$
of class $C^\infty$ relative to the action of $\IT^2$ we define
its left twist $l(T)$ by
\begin{equation}
l(T) = \sum_{n_1,n_2} \, \widehat{T}_{n_1,n_2} \, \lambda^{n_2 p_1} \, ,
\end{equation}
and its right twist $r(T)$ by
\begin{equation}
r(T) = \sum_{n_1,n_2} \,  \lambda^{n_1 p_2} \,  \widehat{T}_{n_1,n_2} \, ,
\end{equation}
Since $|\lambda | = 1$ and $p_1$, $p_2$ are self-adjoint, both series
converge in norm. The construction involves in the case of half-integral spin
 the choice of a square root of $\lambda$.\\
One has the following,
\begin{lemma}\label{lem1} ~\\
{\rm a)}   
Let $x$ be a homogeneous operator of bidegree $(n_1,n_2)$
and $y$ be
a homogeneous operator of  bidegree $(n'_1,n'_2)$.
Define
\begin{equation}
x * y = \lambda^{n'_1 n_2} \, x y \, ; \label{star}
\end{equation}
then $l(x) l(y) = l(x * y)$. \\
{\rm b)}  Let $x$ and $y$ be homogeneous operators as before. Then,
\begin{equation}
l(x) \, r(y) \, - \,  r(y) \, l(x) = (x \, y \, - y \, x) \,
\lambda^{n'_1(n_2 + n'_2)} \lambda^{n_2 p_1 + n'_1 p_2} \, .
\end{equation}
In particular, $[l(x), r(y)] = 0$ if $[x, y] = 0$.
\end{lemma}
\bigskip
\noindent
To check a) and b) one simply uses the following commutation rule which
follows from \eqn{t2act} and it is fulfilled for any
homogeneous operator $T$ of bidegree $(m,n)$,
\begin{equation}
\lambda^{a p_1 + b p_2} \, T = \lambda^{a m + b n} \, T \, \lambda^{a p_1 +
b p_2} \, ,
\, \,\, \forall a, b \in \IZ \, .
\end{equation}
The product $*$ defined in equation (\ref{star}) extends by linearity
to an associative $*$-product on the linear space of smooth operators. \\
One could also define a deformed `right product'. If $x$ is homogeneous of
bidegree
$(n_1,n_2)$ and $y$ is homogeneous of bidegree $(n'_1,n'_2)$ the product is
defined by
\begin{equation}
x *_{r} y = \lambda^{- n'_1 n_2} \, x y \, .
\end{equation}
Then, as with the previous lemma one shows that $r(x) r(y) = r(x
*_{r} y)$.

\bigskip

By Lemma~{\ref{lem1}}~a) one has
that $l(C^{\infty}(S^4))$ is still an algebra and we shall identify it with (the image on
the Hilbert space $\ch$ of) the algebra $C^{\infty}(S_{\theta}^4)$ of smooth functions on the 
deformed sphere $S_{\theta}^4)$. \\
We can then define a new spectral triple $\left(l(C^{\infty}(S^4))\simeq C^{\infty}(S_{\theta}^4), 
\ch , D\right)$ where both the Hilbert space $\ch$ and the operator
$D$ are unchanged while the algebra $C^{\infty}(S^4)$ is modified to $l(C^{\infty}(S^4))\simeq
C^{\infty}(S_{\theta}^4)$. Since
$D$ is of bidegree
$(0,0)$ one has that
\begin{equation}
[D, \,  l(a) ] = l([D, \, a]) \label{bound}
\end{equation}
which is enough to check that $[D, x]$ is bounded for any $x \in  l(\ca)$.

\bigskip
Next, we also deform the real structure by 
twisting the charge conjugation isometry $J$ by
\begin{equation}
\wt{J} = J \, \lambda^{- p_1 p_2} \, . \label{realtwist}
\end{equation}
Due to the antilinearity of $J$ one has that $\wt{J} = \lambda^{p_1 p_2} \, J $ and hence
\begin{equation}
\wt{J}^2  = J ^2 \, .
\end{equation}

\begin{lemma}\label{lem2}~\\
For $x$ homogeneous of bidegree $(n_1,n_2)$ one has that
\begin{equation}
\wt{J} \, l(x) \, \wt{J}^{-1} = r(J \, x \, J^{-1}) \, .
\end{equation}
\end{lemma}
For the proof one needs to check that $\wt{J} \, l(x) = r(J \, x \, J^{-1})  \,  \wt{J}$.
One has
\begin{equation}
\lambda^{-p_1 p_2} \, x = x \, \lambda^{-(p_1 + n_1) (p_2 + n_2)} = x \,
\lambda^{-n_1 n_2} \, \,
\lambda^{-(p_1 n_2 + n_1 p_2)} \, \,  \lambda^{-p_1 p_2}.
\end{equation}
Then
\begin{equation}
\wt{J} \,
l(x) = J \, \lambda^{-p_1 p_2}
\, x \, \lambda^{n_2 p_1} = J
\, x \, \lambda^{-n_1 n_2} \,
\lambda^{- n_1 p_2} \, \, \lambda^{-p_1 p_2},
\end{equation}
 while
\begin{equation}
r(J \, x \, J^{-1}) \, \wt{J} = \lambda^{-n_1 p_2} \, J \, x \,
J^{-1} \, 
J \, \lambda^{-p_1 p_2} = J \, x \, \lambda^{-n_1 ( p_2+n_2)}  \,\, \lambda^{-p_1
p_2}.
\end{equation}
Thus one gets the required equality of Lemma~\ref{lem2}. 

\bigskip
For $x,y \in l(\ca)$ one checks that
\begin{equation}
[x, \, y^0] = 0 \, , \, \, \, y^0 = \wt{J} \, y^* \, \wt{J}^{-1} \, .
\label{commut}
\end{equation}
Indeed, one can assume that $x$ and $y$ are homogeneous and use
Lemma~{\ref{lem2}} together with Lemma~{\ref{lem1}}~a).
Combining equation (\ref{commut}) with equation (\ref{bound}) one then
checks the order one condition
\begin{equation}
[ \, [D, \, x ] \, , y^0 ] = 0 \, , \, \, \, \fl \, x, y \in l(\ca) \,
.\label{first}
\end{equation}
Summing up, we have the following
\begin{theorem}\label{th3}~\\
{\rm a)} The spectral triple $(C^{\infty} (S^{4}_{\theta}),  \ch , D)$
fulfills all axioms of noncommutative manifolds. \\
{\rm b)} Let $e \in C^{\infty} (S^{4}_{\theta} , \Mat_4 (\IC))$ be the
canonical idempotent given in (\ref{idempot}).
The Dirac operator $D$ fulfills
$$
\left\langle \left( e - \frac{1}{2} \right) [D,e]^{4} \right\rangle = \gamma
$$
where $\langle \ \rangle$ is the projection on the commutant of
$\Mat_4 (\IC)$ (i.e. a partial trace) and $\gamma$ is the grading operator.
\end{theorem}
Moreover, the real structure is given by
the twisted involution $\wt{J}$ defined in (\ref{realtwist}). One checks using
the results of \cite{ri93bis} and \cite{co96} that Poincar\'e duality 
continues to hold for the deformed spectral triple. \\
Theorem~\ref{th3} can be extended 
to all metrics on the sphere $S^4$ which are invariant under rotation
of $u$ and $v$ and have the same volume form as the round metric.
In fact, by paralleling the construction for the sphere described above, one can extend it quite
generally
\cite{cola}:
\begin{theorem}\label{th4}~\\
Let $M$ be a compact spin Riemannian manifold whose isometry group has rank
$\geq 2$ (so that one has an inclusion $\IT^2 \subset {\rm Isom}(M)$). Then $M$ admits a
natural one-parameter isospectral deformation to noncommutative (spin) geometries
$M_\theta$.
\end{theorem}
Let $(\ca , \ch , D)$ be the canonical spectral triple associated with a
compact Riemannian spin manifold $M$ as described in Ex.~\ref{ex:ctm}. Here $\ca = \cinf(M)$ is the
algebra of smooth functions on $M$; $\ch= L^2(M,\cs)$ is the Hilbert space of spinors and $D$
is the Dirac operator. Finally, there is the charge conjugation operator $J$, an antilinear
isometry of $\ch$ which gives the real structure. \\
The deformed spectral triple is given by $(l(\ca) , \ch , D)$ with $\ch =
L^2(M,\cs)$ the Hilbert space of spinors, $D$ the Dirac operator and $l(\ca)$ is
really the algebra of smooth functions on $M$ with product deformed to a
$*$-product defined in a way exactly similar to (\ref{star}). The real structure is given by
the twisted involution $\wt{J}$ defined as in (\ref{realtwist}). And again, by 
the results of \cite{ri93bis} and \cite{co96}, Poincar\'e duality 
continues to hold for the deformed spectral triple. 

\subsection{Noncommutative spherical manifolds}\label{ncspman}

As we have seen, on the described deformations one changes the algebra and
the way it acts on the Hilbert
space while keeping the latter and the Dirac operator unchanged, thus getting isospectral 
deformations. From the  decomposition \eqn{t2act} and the deformed product \eqn{star}  one
sees that a central role is played by tori and their noncommutative generalizations. We are
now going to describe in more details this use of the noncommutative tori. 

Let $\theta = (\theta_{jk} = -\theta_{kj})$ be a real antisymmetric $n
\times n$ matrix. 
The noncommutative
torus $\IT_\theta^n$ of `dimension' $n$ and twist $\theta$  is
the `quantum space' whose algebra of polynomial functions 
$A(\IT_\theta^n)$ is
 generated by $n$ independent unitaries
$\row{u}{1}{n}$, subject  to the commutation relations \cite{co80,ri81}
\begin{equation}
u_j u_k = e^{2\pi \ii\theta_{jk}} u_k u_j.
\label{lnct}
\end{equation}
The corresponding $\Cs$-algebra of continuous functions is  
the universal $\Cs$-algebra $C(\IT_\theta^n)$ with the same generators and
relations. 
There is an action $\tau$ of $\IT^n$ on this
$\Cs$-algebra. If $\alpha=(\row{\alpha}{1}{n})\in\IT^n$, this action is given by
$$
\tau(e^{2\pi \ii \alpha}) : u_j \mapsto e^{2\pi \ii \alpha_j} \,u_j.
$$
The smooth subalgebra $\cinf(\IT_\theta^n)$ of $C(\IT_\theta^n)$ under this
action  
consists of
rapidly convergent Fourier series of the form $\sum_{r\in\IZ^n} a_r u^r$,
with $a_r \in \IC$, where
$$
u^r := e^{-\pi \ii r_j\theta_{jk}r_k} \,u_1^{r_1} u_2^{r_2} \dots u_n^{r_n}.
$$
The unitary elements $\{u^r : r \in \IZ^n \}$ form a Weyl system
\cite{grvafi},
since
$$
u^r\,u^s = e^{\pi \ii r_j\theta_{jk}s_k}\, u^{r+s}.
$$
The phase factors 
\begin{equation}
\rho_\theta(r,s) := \exp\{\pi \ii r_j\theta_{jk}s_k\}
\label{phafac}
\end{equation}
form a $2$-cocycle for the group~$\IZ^n$, which is skew (i.e.,
$\rho_\theta(r,r) = 1$) since $\theta$ is skew-symmetric. This also means that  
$C(\IT_\theta^n)$ may be defined as the twisted group $\Cs$-algebra
$C(\IZ^n,\rho_\theta)$.

Let now $M$ be a compact manifold (with no boundary) carrying a smooth action $\sigma$ of
a torus $\IT^n$ of dimension $n \geq 2$. By averaging the translates of
a given Riemannian metric on~$M$ over this torus, we may assume that
$M$ has a $\IT^n$-invariant metric~$g$, so that $\IT^n$ acts by
isometries.

The general $\theta$-deformation of $M$ can be accomplished in two
equivalent ways. The general $\theta$-deformation of $M$ can be accomplished in two
equivalent ways. Firstly, in \cite{cola} the deformation is given
 by a star product of ordinary functions along the lines of \cite{ri93} 
(see also \cite{va}). Indeed, the algebra $\cinf(M)$ may be decomposed
into spectral subspaces which are indexed by the dual group $\IZ^n = \wh\IT^n$. Now, 
each $r \in \IZ^n$ labels a character
$e^{2\pi i\alpha} \mapsto e^{2\pi i r\cdot\alpha}$ of $\IT^n$, with the scalar product
$r\cdot\alpha := r_1\alpha_1 +\cdots+ r_n\alpha_n$. The $r$-th spectral
subspace
for the action $\sigma$ of $\IT^n$ on $\cinf(M)$ consists of those smooth
functions $f_r$ for which
$$
\sigma(e^{2\pi i\alpha}) f_r = e^{2\pi \ii r\cdot\alpha} \,f_r ,
$$
and each $f \in \cinf(M)$ is the sum of a unique (rapidly convergent)
series $f = \sum_{r\in\IZ^n} f_r$.

The $\theta$-deformation of
$\cinf(M)$ may be defined by replacing the ordinary product by a
\textit{Moyal product}, defined on spectral subspaces by
\begin{equation}
f_r \star_\theta g_s := \rho_\theta(r,s)\, f_r g_s ,
\label{moyal-rs}
\end{equation}
with $\rho_\theta(r,s)$ the phase factor in \eqn{phafac}. 
Thus the deformed product is also taken
to respects the $\IZ^n$-grading of functions.  \\
In particular, when $M = \IT^n$ with the obvious
translation action, the algebras $(\cinf(\IT^n),\star_\theta)$ and
$\cinf(\IT_\theta^n)$ are isomorphic.

In the general case, we write $\cinf(M_\theta) := (\cinf(M),\star_\theta)$. Thus, at the
level of smooth algebras the deformation is given explicitly by the star product of ordinary
smooth functions. It is shown in \cite{ri93} that there
is a natural completion of the algebra $\cinf(M_\theta)$ to a $\Cs$-algebra $C(M_\theta)$
whose smooth subalgebra (under the extended action of $\IT^n$) is
precisely the algebra $\cinf(M_\theta)$.

An equivalent approach \cite{codu-vi}, is to
define $C(M_\theta)$ as the fixed-point $\Cs$-subalgebra of
$C(M) \ot C(\IT_\theta^n)$ under the action $\sigma  \times \tau^{-1}$ of $\IT^n$
defined by
$$
e^{2\pi i\alpha} \cdot (f \ot a)
  := \sigma(e^{2\pi i\alpha})\,f \ot \tau(e^{-2\pi i\alpha})\,a;
$$
that is,
$$
C(M_\theta) := \bigl( C(M) \ot C(\IT_\theta^n) \bigr)^{\sigma
\times\tau^{-1}}.
$$
The smooth subalgebra is then given by
\begin{equation}
\cinf(M_\theta)
  := \bigl( \cinf(M) \widehat\otimes \cinf(\IT_\theta^n) \bigr)^{\sigma
\times\tau^{-1}},
\label{fixptalg}
\end{equation}
with $\widehat\otimes$ denoting the appropriate (projective) tensor product of
Fr\'echet algebras. This approach has the advantage that $\cinf(M_\theta)$
may be  determined by generators and relations with the algebra
structure specified by the basic commutation relations
\eqn{lnct} \cite{codu-vi}.

\subsection{The $\theta$-deformed planes and spheres in any dimensions}

We shall briefly describe these classes of spaces while referring to
\cite{codu-vi} for more details.

Let $\theta = (\theta_{jk} = -\theta_{kj})$ be a real antisymmetric $n 
\times n$ matrix. And denote $\lambda^{jk} = e^{2\pi \ii \theta_{jk}}$;
then we have that $\lambda^{kj} = (\lambda^{jk})^{-1}$ and
$\lambda^{jj}=1$.

Let $\ca(\IR^{2n}_\theta)$ be the complex unital $*$-algebra generated by
$2n$ elements $(z^j, z^{k *}, ~j,k=1,\dots,n)$ with relations
\begin{equation}\label{r2ntheta}
z^j z^k = \lambda^{jk} z^k z^j ~, 
~~z^{j*} z^{k*} = \lambda^{jk} z^{k*} z^{j*}, 
~~z^{j*} z^{k} = \lambda^{kj} z^{k*} z^{j*}~,
\end{equation}
with $j,k = 1, \dots,n$. The $*$-algebra $\ca(\IR^{2n}_\theta)$ can be
thought of as the algebra of
complex polynomials on the noncommutative $2n$-plane $\IR^{2n}_\theta$ since
it is a deformation of the commutative $*$-algebra $\ca(\IR^{2n})$ of
complex polynomial functions on $\IR^{2n}$ to which it reduces for
$\theta=0$. 
From relations \eqn{r2ntheta}, it follows that the elements 
$z^{j*} z^{j} = z^{j} z^{j*}~, ~j=1,\dots,n$, are in the center of
$\ca(\IR^{2n}_\theta)$. Since $\sum_{j=1}^n z^{j} z^{j*}$ is central as
well,
it makes sense to  
define $\ca(S^{2n-1}_\theta)$ to be the quotient of the $*$-algebra
$\ca(\IR^{2n}_\theta)$ by the ideal generated by 
$\sum_{j=1}^n z^{j} z^{j*} -1$. 
The $*$-algebra $\ca(S^{2n-1}_\theta)$ can be thought of as the algebra of
complex polynomials on the noncommutative $(2n-1)$-sphere $S^{2n-1}_\theta$
since it is a deformation of the commutative $*$-algebra $\ca(S^{2n-1})$ of
complex polynomial functions on the usual sphere $S^{2n-1}$.

Next, one defines $\ca(\IR^{2n+1}_\theta)$ to be the complex unital
$*$-algebra generated by $2n+1$ elements made of $(z^j, z^{j *},
j=1,\dots,n)$ and of
an addition hermitian element $x=x^*$ with relations like \eqn{r2ntheta} and
in addition
\begin{equation}\label{r2n+1theta}
z^j x = x z^j ~, ~~ j = 1,\dots,n ~. 
\end{equation}
The $*$-algebra $\ca(\IR^{2n+1}_\theta)$ is the algebra 
of complex polynomials on the noncommutative $(2n+1)$-plane
$\IR^{2n+1}_\theta$. 

By the very definition, the elements
$z^{j*} z^{j} = z^{j} z^{j*}~, ~j=1,\dots,n$, and $x$ are in the center of
$\ca(\IR^{2n+1}_\theta)$ and so is the element 
$\sum_{j=1}^n z^{j} z^{j*} +x^2$. 
Then one defines $\ca(S^{2n}_\theta)$ to be the quotient of the
$*$-algebra $\ca(\IR^{2n+1}_\theta)$ by the ideal generated by
$\sum_{j=1}^n z^{j} z^{j*} + x^2 -1$.
The $*$-algebra $\ca(S^{2n}_\theta)$ is the algebra of
complex polynomials on the noncommutative $2n$-sphere $S^{2n}_\theta$ and
is a deformation of the commutative $*$-algebra $\ca(S^{2n})$ of
complex polynomial functions on a usual sphere $S^{2n}$.
By construction the sphere $S^{2n}_\theta$ is a suspension of the sphere
$S^{2n-1}_\theta$.

\bigskip
Next, let $\Cliff(\IR^{2n}_\theta)$ be the unital associative $*$-algebra over $\IC$ generated by $2n$
elements $\Gamma^j, \Gamma^{k *}~, ~j,j=1,\dots,n$, with relations
\begin{eqnarray}
&& \Gamma^j \Gamma^{k} + \lambda^{kj} \Gamma^k \Gamma^{j} = 0~, \nn
&& \Gamma^{j *} \Gamma^{k *} + \lambda^{kj} \Gamma^{k *} \Gamma^{j *} = 0~, \nn
&& \Gamma^j \Gamma^{k *} + \lambda^{jk} \Gamma^k \Gamma^{j *} = \delta^{jk} \II~,  
\end{eqnarray}
where $\II$ is the unit of the algebra and $\delta^{jk}$ is the usual flat metric. For $\theta=0$ one
gets the usual Clifford algebra $\Cliff(\IR^{2n}$ of $\IR^{2n}$. The element
$\gamma\in\Cliff(\IR^{2n}_\theta)$ defined by
\begin{equation}
\gamma = \left[\Gamma^{1 *}, \Gamma^{1}\right]\cdot \ldots \cdot \left[\Gamma^{n *}, \Gamma^{n}\right]
\end{equation}
is hermitian, $\gamma=\gamma^*$, satisfies 
\begin{equation}
\gamma^2=\II, ~~~ \gamma \Gamma^j + \Gamma^j \gamma =0~, ~~\gamma\Gamma^{j *}+\Gamma^{j *}\gamma =0~, 
\end{equation}
and determines a $\IZ_2$-grading of $\Cliff(\IR^{2n}_\theta)$, $\Lambda \mapsto \gamma \Lambda \gamma$.
In fact, one shows \cite{codu-vi} that $\Cliff(\IR^{2n}_\theta)$ is isomorphic to the usual Clifford
algebra $\Cliff(\IR^{2n}$ as a $*$-algebra and as a $\IZ_2$-graded algebra. Furthermore, there is a
representation of $\Cliff(\IR^{2n}_\theta)$ for which
$ \gamma = \left(\begin{smallmatrix}
\II & 0\\
0 & \II 
\end{smallmatrix}\right)$ and $\Gamma^j\in\Mat_{2^n}(\IC)$ of the form
\begin{equation}
\Gamma^j = 
\left(\begin{matrix}
0 & \sigma^j\\
\bar{\sigma}^{j *}& 0  
\end{matrix}\right)~, ~~~
\Gamma^{j *} = 
\left(\begin{matrix}
0 & \bar{\sigma}^{j *}\\
\sigma^j & 0  
\end{matrix}\right)~, ~~~ 
\end{equation}
with $\sigma^j$ and $\bar{\sigma}^{j *}$ in $\Mat_{2^{n-1}}(\IC)$. 

\begin{theorem}~\\
{\rm a)} There is a canonical projection $e\in\Mat_{2^{n}}(\ca(S^{2n}_\theta))$ given by
\begin{equation}
e=\frac{1}{2} \left(\II + \sum_{j=1}^n(\Gamma^{j *} z^j  + \Gamma^j z^{j *} + \gamma x)   \right)~,
\label{proj2n}
\end{equation}
where $(z^j,z^{k *},x)$ are the generators of $\ca(S^{2n}_\theta)$. Moreover, one has that
\begin{equation}
\ChC_k (e) = 0~, \qquad 0 \leq k \leq n-1 \, . \label{vaniproj}
\end{equation}
{\rm a)} There is a canonical unitary $u\in\Mat_{2^{n-1}}(\ca(S^{2n-1}_\theta))$ given by
\begin{equation}
u = \sum_{j+1}^n(\bar{\sigma}^{j} z^j  + \sigma^j z^{j *})    ~,
\label{unit2n}
\end{equation}
where $(z^j,z^{k *})$ are the generators of $\ca(S^{2n-1}_\theta)$. Moreover, one has that
\begin{equation}
\ChC_{k+\frac{1}{2}} (e) = 0~, \qquad 0 \leq k \leq n-1 \, . \label{vaniunit}
\end{equation}
\end{theorem}
For a proof we refer to \cite{codu-vi}.\\
The projection $e$ in \eqn{proj2n} and the unitary $u$ in \eqn{unit2n} provide noncommutative
solutions, via constraints \eqn{vaniproj} and \eqn{vaniunit}, for the algebras
$A_{m,r} = C ({\rm Gr}_{m,r})$ and $B_{m,r} = C ({\rm Gr}_{m,r})$ defined in \eqn{eq2.14}. Thus, there
are admissible surjections 
\begin{equation}
A_{2n,2^{n}} \ra C(S_{\theta}^{2n})~, ~~~B_{{2n-1},2^{n-1}} \ra C(S_{\theta}^{2n-1})~
\end{equation}
The projection \eqn{proj2n} generalizes to higher dimensions the projection constructed in \eqn{idempot}
for the four dimensional sphere $S_\theta$.

\subsection{Gauge theories}
From Theorem~\ref{th4} we know that the deformed spheres, and in particular
the even ones, 
$S^{2n}_\theta$, can be endowed with the structure of a noncommutative spin manifold,
via a spectral triple $(\cinf(S^{2n}_\theta), \ch , D)$, which is isospectral since both the Hilbert
space
$\ch= L^2(S^{2n},\cs)$ and the Dirac operator $D$ are the usual one on the commutative 
sphere $S^{2n}$ whereas only the algebra and its representation on $\ch$ are changed.
In particular one could take the Dirac operator of the usual `round' metric. Then out of this one can
define a suitable Hodge operator $*_H$ on $S^{2n}_\theta$. It turns out that the canonical projection
\eqn{proj2n} satisfies self-duality 
equations
\begin{equation}
*_H e (d e)^n = \ii^n e (d e)^n .
\end{equation}
These equations were somehow `in the air'. For the four dimensional case I
mentioned them during a talk in Ancona in February 2001 \cite{la01}. For
the general case they were derived in \cite{codu-vi} and also in
\cite{asbo}. Their `commutative' counterparts were proposed in
\cite{du-vige,du-vi} together with a description of gauge theories in
terms of projectors.

In particular on the four dimensional sphere $S_\theta$, one can develop
Yang-Mills theory,  since there are all the
required structures, namely the algebra, the calculus (by means of the Dirac operator) and the ``vector
bundle" $e$. \\
The Yang-Mills action is given by, 
\begin{equation}
YM(  \nabla) ={\int \!\!\!\!\!\! -} \,  \theta ^2 \,  ds^{4}~, 
\end{equation}
where $\theta = \nabla^2$ is the curvature, and $ds=D^{-1}$. This action has a strictly positive lower
bound \cite{co94} given by a topological invariant which is just the index \eqn{eq8}, 
\begin{equation}
\varphi(e) = {\int \!\!\!\!\!\! -} \, \gamma (e -\frac{1}{2}) [D, e]^4
\ ds^4 ~. 
\end{equation}
For the canonical projection \eqn{idempot}, owing to \eqn{eq11}, this topological invariant turns out
to be just $1$,
\begin{equation}
\varphi(e)=1~. 
\end{equation}
An important problem, which is still open, is the construction and the classification of Yang-Mills
connections in the noncommutative situation along the line of the ADHM construction \cite{at79}.
For the noncommutative torus this was done in
\cite{cori} and for a noncommutative $\IR^4$ in \cite{nesc}.

\section{Euclidean and Unitary Quantum Spheres} \label{se:qes}
The contents of this Section is essentially a subset of the paper \cite{hala} by Eli 
Hawkins and myself.

The quantum Euclidean spheres in any dimensions, $S_q^{N-1}$, are
(quantum) homogeneous spaces of quantum orthogonal groups, $\SO_q(N)$
\cite{FadResTak89}.  The natural coaction of $\SO_q(N)$ on $\IR^N_q$,
\begin{equation}\label{coactorth} 
\delta : A(\IR^N_q) \rightarrow
A(\SO_q(N)) \otimes A(\IR^N_q), 
\end{equation} 
preserves the `radius of the sphere' and yields a coaction of the quantum 
group $\SO_q(N)$ on $S^{N-1}_q$. \\
Similarly, `odd dimensional' quantum
spheres $S^{2n-1}_q$ can be constructed as  
noncommutative homogeneous spaces of  quantum unitary
groups $\SU_q(n)$ \cite{VakSoi91} (see also \cite{Wel98}). 
Then, analogously to \eqref{coactorth},
there is also a coaction
of the quantum group $\SU_q(n)$ on $S^{2n-1}_q$ 
\begin{equation}\label{coactunit}
\delta : A(S^{2n-1}_q) \rightarrow A(\SU_q(n)) \otimes A(S^{2n-1}_q) .
\end{equation}
In fact, it was realized in \cite{hala} that odd quantum Euclidean spheres
are the same as unitary ones. This fact extends the classical result that
odd dimensional spheres are simultaneously homogeneous spaces of
orthogonal and of unitary groups.

The $*$-algebra $A(S^{N-1}_q)$ of polynomial functions on each of
the spheres $S^{N-1}_q$ 
is given by generators and relations which were expressed in 
terms of a self-adjoint, unipotent matrix 
(a matrix of functions whose square is the identity) which is defined
recursively. Instead in \cite{LanMad01} the algebra was described by means 
of a suitable self-adjoint idempotent
(a matrix of functions whose square is itself).
Let us then describe the algebra $A(S^{N-1}_q)$.
It is generated by elements $\{x_0=x_0*, x_i, x_i^*,  
i=1,\dots,n \}$ for $N=2n+1$ while for $N=2n$ there is no $x_0$.
These generators obey the following commutation relations, 
\begin{subequations}
\label{core}
\begin{align}
\begin{split}
   x_i x_j &= q x_j x_i ,  \quad 0\leq i<j \leq n, \\
   \qquad x_i^* x_j &= q x_j x_i^*, \quad i\neq j ,
\end{split}
\label{core1} \\
   [x_i, x_i^*] &= (1 - q^{-2}) s_{i-1},
\label{core2}
\end{align}
\end{subequations}
with the understanding that $x_0=0$ if $N=2n$, so that in this case the
generator
$x_1$  is normal,
\begin{equation}
x_1x_1^*=x_1^*x_1 \qquad \mathrm{in} \qquad A(S^{2n-1}_q) .
\end{equation}
The `partial radii' $s_i\in A(S^{2n}_q)$, are 
given recursively by
\begin{equation}
\label{pa.ra}
\begin{split}
s_i &:= s_{i-1} + x_i^* x_i = q^{-2} s_{i-1} + x_i x_i^* ,  \\
s_0 &:=x_0^2,
\end{split}
\end{equation}
and the last one $s_n$ which can be shown to be central, is normalized to 
\begin{equation} 
s_n=1
\end{equation}
We see that the equality of the
two formul{\ae} for the elements $s_i$ in \eqref{pa.ra} is equivalent to 
the commutation relation \eqref{core2}.
These $s_i$ are self-adjoint and related as
\begin{equation}
\label{order}
0\leq s_0 \leq \dots \leq s_{n-1} \leq s_n = 1.
\end{equation}
From the commutation relations \eqref{core1} it follows for $i<j$ that
$x_i^*x_ix_j = q^2 x_j x_i^*x_i$; on the other hand $x_j^*x_j x_i =
x_i x_j^*x_j$. By induction, we deduce that
\[
s_i x_j =
\begin{cases} q^2 x_j s_i &: i<j \\
\; x_j s_i &: i\geq j,
\end{cases}
\qquad s_i x_j^* =
\begin{cases} q^{-2} x_j^* s_i &: i<j \\
\; x_j^* s_i &: i\geq j
\end{cases}
\]
and that the $s_i$'s are mutually commuting. They can be used to 
construct representations of the algebra as we shall show later on.

\bigskip
As we have mentioned, in  \cite{LanMad01} it was shown that the defining 
relations of the algebra $A(S^{N-1}_q)$ are equivalent to the condition
that a certain matrix over $A(S^{N-1}_q)$ be idempotent. 
In \cite{hala} it was proven that this  is
also equivalent to the condition that another matrix be unipotent, as we 
shall explain presently.
First consider the even spheres $S^{2n}_q$ for any integer $n>0$. The
algebra $A(S^{2n}_q)$  is generated by elements $\{ x_0, x_i, x_i^*,
i=1, \dots, n\}$. Let us first consider the free unital $*$-algebra $F:=\IC\langle 1, x_0,
x_i, x_i^*, i=1, \dots, n\rangle$ on $2n+1$ generators.
We recursively define self-adjoint matrices $u_{(2n)}\in
\Mat_{2^n}(F)$ for all
$n$ by,
\begin{equation}
\label{uni.even} u_{(2n)} :=
\begin{pmatrix} q^{-1}u_{(2n-2)} & x_n \\ x_n^* & -u_{(2n-2)}
\end{pmatrix}
,\end{equation} with $u_{(0)}=x_0$. The $*$-algebra $A(S^{2n}_q)$
is then defined by the relations that $u_{(2n)}$ is unipotent,
$u_{(2n)}^2=1$, and self-adjoint, $u_{(2n)}^*=u_{(2n)}$. That is, the algebra is the quotient of 
the free algebra $F$ by these relations.\\
The self-adjointness relations merely give that $x_i^*$ is the adjoint
of $x_i$ and $x_0$ is self-adjoint. Unipotency gives a
matrix of $2^{2n}$ relations, although many of these are vacuous or
redundant. These can be deduced inductively from \eqref{uni.even}
which gives,
\begin{equation}
\label{u2recursion}
   u_{(2n)}^2 =
\begin{pmatrix} q^{-2}u_{(2n-2)}^2 + x_nx_n^* & q^{-1}u_{(2n-2)} x_n -
x_n u_{(2n-2)} \\ q^{-1}x_n^* u_{(2n-2)} - u_{(2n-2)} x_n^* &
u_{(2n-2)}^2 + x_n^*x_n
\end{pmatrix}.
\end{equation}
The condition that $u_{(2n)}^2=1$ means in particular that
$u_{(2n)}^2$ is
diagonal with all the diagonal entries equal. Looking at
\eqref{u2recursion}, we
see that the same must be true of $u_{2n-2}^2 \in
\Mat_{2^{n-1}}(A(S^{2n}_q))$,
and so on. Thus, the diagonal relations require that all the diagonal
entries of
(each) $u_{(2j)}^2$ are equal. If this is true for $u_{(2j-2)}^2$,
then the relation
for $u_{(2j)}^2$ is that the same element (the diagonal entry) can be
written in
two different ways. This element is simply $s_j$ and the two ways of
writing it are
those  given in \eqref{pa.ra}. 
Finally, $u_{(2n)}^2=1$ gives the relation $s_n=1$.
The off-diagonal relations are $q^{-1}u_{(2j-2)}x_j = x_ju_{(2j-2)}$
and
$q^{-1}x_j^*u_{(2j-2)} = u_{(2j-2)}x_j^*$ for every $j = 1, \dots ,
n$. Because
the matrix $u_{(2j-2)}$ is constructed linearly from all of the
generators $x_i$
and $x_i^*$ for $i<j$, this conditions are equivalent to the
commutation relations
\eqref{core1}. 
This presentation of the relations by the unipotency of $u_{(2n)}$ is
also the easiest way to see that there is an
isomorphism $A(S^{2n}_{1/q}) \cong A(S^{2n}_q)$ which is obtained by the 
substitutions $q \leftrightarrow q^{-1}$, $x_0 \rightarrow (-q)^{n} x_0$, 
and $x_i \rightarrow (-q)^{n-i} x_i^*$. This transforms the  
matrix $u_{(2n)} \rightarrow \wt{u}_{(2n)}$ and the latter is 
is unipotent and self-adjoint if and only if $u_{(2n)}$ is. Thus there is 
an isomorphism  $A(S^{2n}_{1/q}) \cong A(S^{2n}_q)$,
and we can assume that $\abs{q}>1$ without loss of generality.

Next, consider the odd spheres $S^{2n-1}_q$ for any integer $n>0$. We
can construct a unipotent $u_{(2n-1)} \in
\Mat_{2^{n}}[A(S^{2n-1}_q)]$, simply by setting
$x_0=0$ in $u_{(2n)}$. Once again, the unipotency condition,
$u_{(2n-1)}^2=1$, is equivalent to the relations defining the algebra
$A(S^{2n-1}_q)$ of polynomial functions on $S^{2n-1}_q$. Again, one
defines self-adjoint elements $s_i\in A(S^{2n-1}_q)$ such  that $s_i =
s_{i-1} + x_i^*x_i = q^{-2} s_{i-1} + x_i x_i^*$ with now $s_0 =
x_0^2 =0$. The commutation relations are again given by \eqref{core} but 
now \eqref{core2} gives in particular that the generator
$x_1$ is normal, $x_1x_1^*=x_1^*x_1$. 
The previous argument also shows that $A(S^{2n-1}_q)$ is the quotient of 
$A(S^{2n}_q)$ by the ideal generated by
$x_0$. 
Geometrically, we may think of  $S^{2n-1}_q$ as a noncommutative 
subspace of $S^{2n}_q$. Because of the  isomorphism $A(S^{2n}_{1/q})
\cong  A(S^{2n}_q)$, we have another isomorphism $A(S^{2n-1}_{1/q}) \cong 
A(S^{2n-1}_q)$, and again we can  assume that $\abs{q}>1$ without any 
loss of generality.

\begin{remark}
The algebras of our spheres, both in even and odd `dimensions', are generated by the entries
of a projections. This is the same as the condition of {\it full projection} used by S.
Waldmann in his analysis of Morita equivalence of star products \cite{wa}.
\end{remark}

There is also a way of realizing even spheres as noncommutative
subspaces of odd ones. Consider $S^{2n+1}_q$, set $x_1=x_1^*=x_0$ and
relabel $x_2$ as $x_1$, \emph{et cetera}; let $u'_{(2n+1)}$ be the matrix
obtained from $u_{(2n+1)}$ with these substitutions. The matrix
$u'_{(2n+1)}$ is the  same as $u_{(2n)}$ in which we substitute
\[
x_0 \rightarrow {\textstyle \begin{pmatrix}0 &
x_0\\ x_0& 0 \end{pmatrix}},
\qquad x_j \rightarrow
{\textstyle\begin{pmatrix} x_{j}
& 0 \\ 0 & x_{j} \end{pmatrix}}, \quad j\not= 0 \, .
\]
Then the unipotency of $u'_{(2n+1)}$ yields precisely
the same relations coming from the unipotency of $u_{(2n)}$. This
shows
that $A(S^{2n}_q)$ is
the quotient of $A(S^{2n+1}_q)$ by the $*$-ideal generated by
$x_1-x_1^*$. Geometrically, we may think of  $S^{2n}_q$ as a noncommutative 
subspace of$S^{2n+1}_q$.

Summing up, every sphere contains a smaller sphere of dimension one less;
by following this tower of inclusions to its  base, we see that every
sphere contains a classical $S^1$, because the circle does not deform.
From this, it is easy to see that the spheres $S^{N-1}_q$ have a $S^1$ 
worth of classical points. Indeed, with $\lambda \in \IC$ such that $|\lambda|^2=1$, there is a 
family 
of $1$-dimensional representations (characters) of the algebra
$A(S^{N-1}_q)$ given by
\begin{eqnarray}\label{1reps}
&& \psi_{\lambda}(1) = 1, 
~~~~~\psi_{\lambda}(x_{n}) =  \lambda ,
~~~~~\psi_{\lambda}((x_{n})^*) = \bar{\lambda}~, \nn
&& \psi_{\lambda}(x_{i}) = \psi_{\lambda}((x_{i})^*) = 0 ~, 
\end{eqnarray}
for $i=0,1,\ldots, n-1$ or $i=1,\ldots, n-1$ according to whether $N=2n+1$ 
or $N=2n$, respectively. 

Each even sphere algebra has  an involutive
automorphism
\begin{equation}\label{sigma}
\begin{split}
\sigma :\, & A(S^{2n}_q) \to A(S^{2n}_q) \\
& x_0 \mapsto -x_0 ; \qquad x_j \mapsto x_j,
\quad  j\neq 0 ,
\end{split}
\end{equation}
which corresponds to flipping (reflecting) the classical
$S^{2n}$ across the hyperplane $x_0=0$. The coinvariant algebra of
$\sigma$ is the quotient of $A(S^{2n}_q)$ by the ideal generated by
$x_0$, which, as we have noted, is simply
$A(S^{2n-1}_q)$. \\
Geometrically this means that $S^{2n-1}_q$ is the ``equator'' of 
$S^{2n}_q$, the subspace fixed by the flip.

As for odd spheres, they have an action $\rho:\IT\to\Aut[A(S^{2n})]$
of the torus group $\IT$, defined by
multiplying $x_1$ by a phase
and leaving the  other generators unchanged,
\begin{equation}\label{rho}
\begin{split}
\rho(\lambda) :\, & A(S^{2n+1}_q) \to A(S^{2n+1}_q) \\
& x_1 \mapsto \lambda x_1 ; \qquad x_j \mapsto x_j,
\quad  j\neq 1 .
\end{split}
\end{equation}
The coinvariant algebra is
given by setting $x_1=0$. Now, let  $u''_{(2n+1)}$ be the matrix
obtained by
setting $x_1=0$ and relabeling $x_2$ as $x_1$,
\emph{et cetera}, in the matrix  $u_{(2n+1)}$. Then, $u''_{(2n+1)}$ is
equivalent to
tensoring $u_{(2n-1)}$ with  $\left(\begin{smallmatrix}1 &
0\\ 0&1\end{smallmatrix}\right)$,
\[
u''_{(2n+1)} = u_{(2n-1)} \otimes
\left(\begin{smallmatrix}1 &
0\\ 0&1\end{smallmatrix}\right)
\]
and the result is unipotent if and only if
$u_{(2n-1)}$ is; that is the unipotency of $u''_{(2n+1)}$ yields all
and only
the same relations coming from the unipotency of $u_{(2n-1)}$. This
shows
that
$A(S^{2n-1}_q)$ is the quotient of
$A(S^{2n+1}_q)$ by
the $*$-ideal generated by $x_1$ and $S^{2n-1}_q$ is
the noncommutative
subspace of $S^{2n+1}_q$ fixed by the $\IT$-action in \eqref{rho}.

\subsection{The structure of the deformations}
For each deformed sphere  $S^{N-1}_q$, we have a one
parameter family of
algebras $A(S^{N-1}_q)$ which, at $q=1$, gives
$A(S^{N-1}_1)=A(S^{N-1})$, the algebra of
polynomial functions on a classical sphere $S^{N-1}$. It is possible
to identify this one-parameter family of algebras to a fixed vector
space and view the product as varying with the parameter: let us indicate
this product with the symbol 
$*_q$. We can then construct a Poisson bracket on $A(S^{N-1})$ from
the first derivative of the product at the classical parameter value,
$q=1$,
\begin{equation}\label{pobra}
\{f,g\} := -i \left.\frac{d}{dq} ~\right|_{q=1} (f *_q g - g *_q f)~.
\end{equation}
The usual properties of a Poisson bracket (Leibniz and Jacobi
identities) are simple consequences of associativity. 

In general, given such a one-parameter deformation from a commutative
manifold $\cm$ into noncommutative algebras, we can construct a
Poisson bracket on functions.
This Poisson algebra, $A(\cm)$ with the commutative product and the
Poisson bracket, describes the deformation to first order.
A deformation is essentially a path through an enormous space of
possible algebras, and the Poisson algebra is just a tangent.
Nevertheless, if the deformation is well behaved the Poisson algebra
does indicate where it is heading. 

We shall also use the fact that the manifold $\IR^{2n}$ has a unique
symplectic structure, modulo isomorphism. Then, this symplectic structure
corresponds to an essentially unique deformation. If we complete to a
\Cs-algebra, then the deformation
of $\cc_0(\IR^{2n})$ (continuous functions vanishing at infinity) will be
the algebra,
$\ck$, of compact operators on a
countably infinite-dimensional Hilbert space.

Let us go back to the spheres $S^{N-1}_q$ and look more closely at them.

We have seen that the $S^{2n-1}_q$ noncommutative subspace of
$S^{2n}_q$
corresponds to the equator, $S^{2n-1}\subset S^{2n}$, where
$x_0=0$ and the
Poisson structure  on $S^{2n}$ is degenerate. On the remaining
$S^{2n}\smallsetminus S^{2n-1}$, the Poisson structure is
nondegenerate.
So, topologically, we have a union of two copies of symplectic
$\IR^{2n}$. Then, the kernel of the quotient map $A(S^{2n}_q)\to
A(S^{2n-1}_q)$ should be a deformation of the subalgebra of functions on
$S^{2n}$ which
vanish at the equator. If we complete to \Cs-algebras, this should
give us the  direct sum of two copies of $\ck$, one for each
hemisphere.
Thus we expect that the
\Cs-algebra $\cc(S^{2n}_q)$ will be an extension:
\begin{equation}\label{ext.even}
0 \to \ck\oplus\ck \to \cc(S^{2n}_q) \to \cc(S^{2n-1}_q)
\to 0
.\end{equation}

In odd dimensions, the Poisson structure is necessarily degenerate.
However, the
$S^{2n-1}_q$ noncommutative subspace of $S^{2n+1}_q$ corresponds
classically to the Poisson structure being more degenerate on
$S^{2n-1}\subset  S^{2n+1}$. It is of rank $2n$ at most points, but of
rank $2n-2$ (or less) along
$S^{2n-1}$. The complement
$S^{2n+1}\smallsetminus S^{2n-1}$ has a symplectic foliation by
$2n$ dimensional leaves which is invariant under the $\IT$ action;
the simplest possibility is that this corresponds to the product in
the
identification
\[
S^{2n+1}\smallsetminus S^{2n-1} \cong S^1\times \IR^{2n}
.\]
If we complete to \Cs-algebras, then the deformation of this should give the algebra  
$\cc(S^1)\otimes\ck$. The kernel of the quotient map $A(S^{2n+1}_q)\to
A(S^{2n-1}_q)$ should be this deformation, so we expect another
extension,
\begin{equation}
\label{ext.odd} 0 \to \cc(S^1)\otimes\ck \to \cc(S^{2n+1}_q)\to
\cc(S^{2n-1}_q) \to 0
.\end{equation}

The extensions \eqn{ext.even} and \eqn{ext.odd} turn out to be correct. As we have
mentioned, the odd dimensional spheres we are considering are equivalent
to the ``unitary" odd quantum spheres of Vaksman and Soibelman
\cite{VakSoi91}. In
\cite{HonSzy01} Hong and Szyma{\'n}ski obtained the \Cs-algebras
$\cc(S^{2n+1}_q)$ as Cuntz-Krieger algebras of suitable graphs. From
this
construction they derived the extension \eqref{ext.odd}. They
also considered even spheres, defined as quotients of odd ones  by
the ideal
generated by $x_1 - x_1^*$. These are thus isomorphic to the even
spheres we are
considering here. They also obtained these as Cuntz-Krieger algebras
and derived
the extension \eqref{ext.even}. However, as explicitly stated in the
introduction to \cite{HonSzy01}, they were unable to realize even
spheres as
quantum homogeneous spaces of quantum orthogonal groups, thus also
failing to
realize that ``unitary" and ``orthogonal" odd quantum spheres are the
same.

\subsection{Representations}

We shall now exhibit
all representations of the algebra $A(S^{N-1}_q)$ which in turn
extend to
the \Cs-algebra $\cc(S^{N-1}_q)$. \\
Representations of the odd dimensional spheres were constructed in
\cite{VakSoi91}. The primitive spectra of all these
spheres were
compute in \cite{HonSzy01},  which amounts to a classification of
representations. The
representations for
quantum Euclidean spheres have also been constructed in \cite{Fio95}
by thinking
of them as quotient algebras of quantum Euclidean planes.

The structure of
the representations can be anticipated from the construction of
$S^{N-1}_q$
via the extensions \eqref{ext.even} and \eqref{ext.odd} and by
remembering that an
irreducible representation $\psi$ can be partially characterized by its
kernel. Moreover, an irreducible representation of a \Cs-algebra
restricts either to an irreducible or a trivial representation of any
ideal; and conversely, an irreducible representation of an ideal
extends to an irreducible representation of the \Cs-algebra (see for
instance \cite{Dav96}).

For an even sphere $S^{2n}_q$, the kernel of an irreducible
representation $\psi$ will
contain one or both of the copies of $\ck\subset \cc(S^{2n}_q)$.
If $\ck\oplus\ck\subseteq\ker\psi$, then $\psi$ factors through
$\cc(S^{2n-1}_q)$ and is given by a representation of that algebra.
If one copy of $\ck$ is not in $\ker\psi$, then $\psi$
restricts to a representation of this $\ck$. However, $\ck$ has only one
irreducible representation. Since $\ck$ is an ideal in $\cc(S^{2n}_q)$,
the unique irreducible representation of
$\ck$ uniquely extends to a representation of $\cc(S^{2n}_q)$ (with the
other copy of $\ck$
in its kernel).
\\
Thus, up to isomorphism the irreducible representations of $S^{2n}_q$ should be: \\
1. all irreducible representations of $S^{2n-1}_q$, \\ 
2. a unique representation with kernel the second
copy of $\ck$, \\
3. a unique representation with kernel the first copy of $\ck$. \\
From the extension \eqref{ext.even} we expect that the generator
$x_0$ is a
self-adjoint  element of $\ck\oplus\ck\subset \cc(S^{2n}_q)$ and it
should have almost
discrete, real spectrum: it will therefore be used to  
decompose the Hilbert space in  a representation.

Similarly, from the construction of $S^{2n+1}_q$ by the extension
\eqref{ext.odd}, one can anticipate the structure of its
representations.
Firstly, if $\cc(S^1)\otimes\ck\subseteq\ker\psi$, then $\psi$ factors
through
$\cc(S^{2n-1}_q)$ and is really a representation of $S^{2n-1}_q$.
Otherwise, $\psi$ restricts to an irreducible representation of
$\cc(S^1)\otimes\ck$. This factorizes as the tensor product of an
irreducible representation of $\cc(S^1)$ with one of $\ck$. The
irreducible representations of $\cc(S^1)$ are simply given by the
points of
$S^1$, and as we have
mentioned, $\ck$ has a unique irreducible representation. The
representations of
$\cc(S^1)\otimes\ck$ are thus classified by the points of $S^1$. These
representations extend  uniquely from the ideal $\cc(S^1)\otimes\ck$ to
the whole algebra $\cc(S^{2n+1}_q)$.
\\
Thus, up to isomorphism, the irreducible representations of $S^{2n+1}_q$ should be:\\
1. all irreducible representations of $S^{2n-1}_q$, \\
2. a family of representations parameterized by $S^1$.

\bigskip
In the construction of the representations, a simple identity regarding  the spectra of operators 
will be especially useful  (see, for instance \cite{Sak98}). If $x$ is an element of any
\Cs-algebra, then
\begin{equation}
\label{spec.id}
\{0\}\cup \Spec x^*x = \{0\}\cup \Spec xx^*
.\end{equation}

\subsection{Even sphere representations} 

To illustrate the general structure we shall start by describing the  lowest
dimensional case, namely $S^2_q$. This is isomorphic to the so-called equator
sphere of Podle{\'s} \cite{Pod87}. For this  sphere, the 
representations were also
constructed in \cite{MasNakWat91} in a way  close to the one presented here.

Let us then consider the sphere $S^2_q$.
\\
As we have discussed, we expect that, in some faithful representation, $x_0$ is a compact operator 
and
thus has an almost discrete, real spectrum. However, we cannot assume
\emph{a priori} that $x_0$ has eigenvalues, let alone that
its eigenvectors
form a complete basis of the Hilbert space. 
The sphere relation $1 = x_0^2 + x_1^*x_1 = q^{-2} x_0^2 + x_1 x_1^* $ shows
that
$x_0^2 \leq 1$ and  thus
$\norm{x_0}
\leq 1$. As $x_0$ is self-adjoint, this shows that $\Spec x_0
\subseteq [-1,1]$. By \eqref{spec.id} we have also,
\begin{align*}
\{0\}\cup \Spec x_1^*x_1 &= \{0\}\cup \Spec x_1x_1^* \\
\{0\}\cup \Spec (1-x_0^2) &= \{0\}\cup \Spec (1- q^{-2}x_0^2) \\
\{1\}\cup \Spec x_0^2 &= \{1\} \cup q^{-2} \Spec x_0^2
.\end{align*}
Because we have assumed that $\abs{q}\geq1$, the only subsets of
$[0,1]$ that
satisfy this condition are $\{0\}$ and
$\{0,q^{-2k}\mid k=0,1,\dots\}$.
\\
If $x_0\neq 0 \in \cc(S^2_q)$ then $\Spec x_0^2$ is the latter set. We
cannot
simply assume that $x_0\neq 0$, since not every $*$-algebra is a
subalgebra of a
\Cs-algebra; however, our explicit representations will show that
that is the
case here.

Now let $\ch$ be a separable Hilbert space and suppose that we have an irreducible
$*$-representation, 
$\psi : A(S^2_q) \to
\cl(\ch)$. \\
If $\psi(x_0)= 0$ then $1 = \psi(x_1)\psi(x_1)^*
=\psi(x_1)^*\psi(x_1)$.
Thus
$\psi(x_1)$ is unitary, and by the assumption of
irreducibility, it is a number $\lambda\in\IC$, $\abs\lambda = 1$. So,
$\ch=\IC$ and the representation is $\psi^{(1)}_\lambda$ defined by,
\begin{equation}\label{reps1}
\psi^{(1)}_\lambda(x_0) =0; \qquad
\psi^{(1)}_\lambda(x_1) =\lambda , \quad \lambda\in S^1 .
\end{equation}
Thus we have an $S^1$ worth of representations with $x_0$ in the
kernel. \\
If $\psi(x_0)\neq 0$, then $1 \in \Spec x_0^2$; it is an isolated point in the spectrum and 
therefore an eigenvalue. For some sign
$\pm$ there exists a unit vector $\ket0 \in\ch$ such that
$\psi(x_0)\ket0 =
\pm \ket 0$. The relation $x_0 x_1=qx_1x_0$
suggests
that $x_1$ and
$x_1^*$ shift the eigenvalues of $x_0$. Indeed, for $k=0,1,\dots$, the vector $\psi(x_1^*)^k\ket 0$ is an
eigenvector as well, because
\[
\psi(x_0)  \psi(x_1^*)^k\ket 0  = q^{-k} \psi({x_1^*}^k x_0) \ket 0  =
\pm q^{-k}
\psi(x_1^*)^k\ket 0
.\] By normalizing, we obtain a sequence of unit
eigenvectors, defined by
\[
\ket k := (1-q^{-2k})^{-1/2} \psi(x_1^*)\ket{k-1}
.\]
We have thus two representations $\psi^{(2)}_+$ and $\psi^{(2)}_-$,
and direct computation shows that
\begin{equation}\label{reps2}
\begin{split}
       \psi^{(2)}_\pm(x_0) \ket k &= \pm
       q^{-k} \ket k,  \\
\psi^{(2)}_\pm(x_1) \ket k &= (1 - q^{-2k})^{1/2}
\ket{k-1} ,  \\
\psi^{(2)}_\pm(x_1^*) \ket k &= (1 - q^{-2(k+1)})^{1/2} \ket{k+1} .
\end{split}
\end{equation}
The eigenvectors $\{\ket k\mid k=0,1,\dots\}$
are mutually orthogonal because they have distinct eigenvalues, and by
the assumption of irreducibility they form a basis for the Hilbert space $\ch$.
\\
Notice that any power of $\psi^{(2)}_\pm(x_0)$ is a trace class
operator, while this is not the case for the operators
$\psi^{(2)}_{\pm}(x_1)$ and $\psi^{(2)}_{\pm}(x_1^*)$ nor for any of
their powers.
\\
Note also that the representations \eqref{reps2} are related by the
automorphism
$\sigma$ in
\eqref{sigma}, as
\begin{equation}
\psi^{(2)}_\pm\circ\sigma = \psi^{(2)}_\mp.
\end{equation}

If we set a value of $q$ with $\abs{q}<1$ in \eqref{reps2}, the
operators would be
unbounded. This is the reason for assuming that $\abs{q}>1$. The
assumption
was used in computing
$\Spec x_0$. Not only is $\norm{x_0}\leq 1$, but by a similar
calculation
$\norm{x_0}\leq \abs{q}$. Which bound is more relevant obviously
depends on
whether $q$ is greater or less than $1$.
For $\abs{q}<1$ the appropriate formul{\ae} for the representations
can be obtained from \eqref{reps2}
by replacing the index $k$ with $-k-1$. As a consequence, the role of
$x_1$
and $x_1^*$ as lowering and raising operators is exchanged.

\bigskip
For the general even spheres $S^{2n}_q$ the structure of the representations is similar to that for
$S^2_q$, but more complicated. The element $x_0$ is no longer sufficient to completely decompose the
Hilbert space of the representation and we need to use all the commuting self-adjoint elements $s_i\in
A(S^{2n}_q)$ defined in \eqref{pa.ra}.  

Suppose that $\psi : A(S^{2n}_q) \to \cl(\ch)$ is an irreducible
$*$-representation. \\
If $\psi(x_0) =0$, then $\psi$ factors through $A(S^{2n-1}_q)$. Thus
$\psi$ is an irreducible representation of $A(S^{2n-1}_q)$; these
will be
discussed later. \\
If $\psi(x_0) \neq 0$, then $\psi(s_0)\neq 0$, and by the relations
\eqref{order}, all of the $\psi(s_i)$'s are nonzero.  Proceeding recursively, we find that
there is a simultaneous eigenspace with
eigenvalue $1$ for all the $\psi(s_i)$'s. That is, there must exist a
unit vector
$\ket{0,\dots,0}\in \ch$ such that $\psi(s_i) \ket{0,\dots,0} =
\ket{0,\dots,0}$ for all
$i$ and $\psi(x_0) \ket{0,\dots,0} = \pm \ket{0,\dots,0}$. More unit
vectors are defined by 
\[
\ket{k_0,\dots,k_{n-1}} \sim \psi(x_1^*)^{k_0}\dots \psi(x_n^*)^{k_{n-1}} \ket{0,\dots,0}
\]
modulo a positive normalizing factor. Working out the correct normalizing factors we get two
representations $\psi^{(2n)}_\pm$ defined by,
\begin{align}
\label{rep.even}
\psi^{(2n)}_\pm(x_0) \ket{k_0,\dots,k_{n-1}} &= \pm
q^{-(k_0+\dots+k_{n-1})}
\ket{k_0,\dots,k_{n-1}}  \\
\psi^{(2n)}_\pm(x_i) \ket{k_0,\dots,k_{n-1}} &= (1-
q^{-2k_{i-1}})^{1/2}
q^{-(\sum_{j=i}^{n-1}k_j)} \ket{\dots,k_{i-1}-1,\dots} \nn
\psi^{(2n)}_\pm(x_i^*) \ket{k_0,\dots,k_{n-1}} &= (1-
q^{-2(k_{i-1}+1)})^{1/2}
q^{-(\sum_{j=i}^{n-1}k_j)} \ket{ \dots,k_{i-1}+1,\dots } \nonumber
\end{align}
with $i=1,\dots,n$. With this values for the index $i$, we see that $x_i$ lowers $k_{i-1}$,
whereas $x_i^*$ raises $k_{i-1}$. From irreducibility  the
collection of
vectors
$\{ \ket{k_0,\dots,k_{n-1}} , k_i \geq 0 \}$ constitute a complete
basis for the Hilbert space $\ch$.  \\
As before, the two representations \eqref{rep.even} are related by
the automorphism $\sigma$, 
\begin{equation}
\psi^{(2n)}_\pm\circ\sigma = \psi^{(2n)}_\mp.
\end{equation}
Again the formul{\ae} \eqref{rep.even} for the representations are corrected for
$\abs{q}>1$; and again the representations for $\abs{q}<1$ can be
obtained by replacing all indices $k_i$ with $-k_i-1$ in \eqref{rep.even}.

In all of the irreducible representations of $A(S^{2n}_q)$, the
representative of $x_0$ is compact; in fact it is trace class. We can
deduce from
this that the \Cs-ideal generated by $\psi^{(2n)}_\pm(x_0)$ in
$\cc(S^{2n}_q)$ is
isomorphic to  $\ck(\ch)$, the ideal of all compact operators on
$\ch$. By using the
continuous functional calculus, we can apply any function
$f\in\cc[-1,1]$ to $x_0$. If $f$ is supported on $[0,1]$, then $f(x_0)
\in \ker \psi^{(2n)}_-$. Likewise if $f$ is supported in $[-1,0]$,
then
$f(x_0) \in \ker \psi^{(2n)}_+$.
     From this we deduce that the \Cs-ideal generated by $x_0$ in
$\cc(S^{2n}_q)$ is $\ck\oplus\ck$. One copy of $\ck$ is $\ker
\psi^{(2n)}_+$; the other is $\ker \psi^{(2n)}_-$. Thus we get exactly  
the extension \eqref{ext.even}.

\subsection{Odd sphere representations}
Again, to illustrate the general strategy we shall work out in detail 
the simplest
case, that of the sphere $S^3_q$. This can be identified with the underlying
noncommutative space of the quantum group $\SU_q(2)$  and as such the
representations of the algebra are well known \cite{wo87}.

The generators $\{x_i, x_i^* \mid
i=1,2 \}$ of the algebra $A(S^3_q)$ satisfy the commutation relations
$x_1 x_2 = q x_2 x_1$, $x_i^* x_j = q x_j x_i^*$,$i\neq j$,
$[x_1, x_1^*] = 0$, and $[x_2, x_2^*] = (1-q^{-2})x_1 x_1^*$.
Furthermore,
there is also the sphere relation $1 = x_2^* x_2 + x_1^* x_1 =  x_2 x_2^* +
q^{-2} x_1 x_1^*$.
\\
The normal generator $x_1$ plays much the same role for the
representations of
$S^3_q$ that $x_0$ does for those of $S^2_q$. The sphere relation
shows that
$\norm{x_1}\leq 1$ and
\begin{align*}
\{0\} \cup \Spec x_2^*x_2 &= \{0\} \cup \Spec x_2x_2^* \\
\{0\} \cup \Spec (1 - x_1^*x_1) &= \{0\} \cup \Spec (1 - q^{-2} x_1
x_1^*) \\
\{1\} \cup \Spec x_1^*x_1 &= \{1\} \cup q^{-2} \Spec x_1^*x_1
,\end{align*}
which shows that either $x_1=0$ or $\Spec x_1^*x_1 =
\{0,q^{-2k}\mid k=0,1,\dots\}$.

Let $\psi : A(S^3_q) \to
\cl(\ch)$ be an irreducible $*$-representation.\\
If $\psi(x_1)=0$ then
the relations reduce to $1 =
\psi(x_2)\psi(x_2)^* =
\psi(x_2)^*\psi(x_2)$. Thus $\psi(x_2)$ is
unitary and by the
assumption of irreducibility, it is a scalar,
$\psi(x_2)=\lambda \in
\IC$ with
$\abs\lambda=1$. Thus, as before, we have an
$S^1$ of representations of this kind.\\
If $\psi(x_1)\neq0$, then
$1\in \Spec \psi(x_1^*x_1)$ and is an isolated point in the spectrum. Thus, there exists
a unit vector $\ket0
\in \ch$ such that $\psi(x_1^*x_1) \ket 0 = \ket
0$, and by the
assumption of irreducibility, there is some
$\lambda\in\IC$
with
$\abs\lambda=1$ such that $\psi(x_1)\ket0 =
\lambda\ket0$. We see then that
$\psi(x_2^*)^k\ket0$ is an
eigenvector
\[
\psi(x_1) \psi(x_2^*)^k\ket0  = q^{-k}
\psi({x_2^*}^kx_1)\ket0  =
\lambda q^{-k}
\psi(x_2^*)^k\ket0
.\]
By normalizing, we get a sequence of unit
eigenvectors
recursively defined by
\[
\ket k := (1-q^{-2k})^{-1/2}
\psi(x_2^*) \ket{k-1}
.\]
A family of
representations $\psi^{(3)}_\lambda$, $\lambda \in S^1$, is then
defined by
\begin{align}\label{reps3}
\psi^{(3)}_{\lambda}(x_1) \ket{k} &=
\lambda q^{-k} \ket{k}
, \nn 
\psi^{(3)}_{\lambda}(x_1^*)
\ket{k} &=  \bar{\lambda} q^{-k} \ket{k}
, \nn 
\psi^{(3)}_{\lambda}(x_2) \ket{k} &=  (1 - q^{-2k})^{1/2}
\ket{k-1}
, \nn 
\psi^{(3)}_{\lambda}(x_2^*) \ket{k} &=  (1 -
q^{-2(k+1)})^{1/2}
\ket{k+1} .
\end{align}
We notice that
any power of $\psi^{(3)}_{\lambda}(x_1)$
or
$\psi^{(3)}_{\lambda}(x_1^*)$ is a trace class operator, while
this is
not the case for the operators $\psi^{(3)}_{\lambda}(x_2)$ and
$\psi^{(3)}_{\lambda}(x_2^*)$ nor for any of their powers.

\bigskip

Next, for  the general odd spheres $S^{2n+1}_q$ , 
let $\psi:A(S^{2n+1}_q)\to\cl(\ch)$ be an irreducible
representation. \\
If $\psi(x_1) =0$ then $\psi$ factors through $A(S^{2n-1}_q)$ and is
an
irreducible representation of that algebra.\\
If $\psi(x_1)\neq 0$ then $\psi(s_1)\neq 0, \psi(s_2)\neq 0, $
\emph{et
cetera}. By the same
arguments as  for $S^{2n}_q$, there must exist a simultaneous
eigenspace with eigenvalue
$1$  for all of $s_1,\dots s_n$. By the assumption of irreducibility,
this
eigenspace is
$1$-dimensional.
Let
$\ket{0,\dots,0} \in \ch$ be a unit vector in this eigenspace. Then
$s_i
\ket{0,\dots,0}
=
\ket{0,\dots,0}$ for $i=1,\dots,n$. The restriction of $\psi(x_1)$ to
this
subspace is unitary and thus for some $\lambda\in\IC$ with
$\abs\lambda=1$, one has that 
$\psi(x_1) \ket{0,\dots,0} = \lambda \ket{0,\dots,0}$.
We can construct more simultaneous eigenvectors of the $s_i$'s by 
defining 
\[
\ket{k_1,\dots,k_n} \sim \psi(x_2)^{k_1}\dots \psi(x_{n+1})^{k_n}
\ket{0,\dots,0}
\]
modulo a positive normalizing constant. 
Working out the normalization, one has a family of representations
$\psi^{(2n+1)}_\lambda$,
\begin{equation}
\begin{split}
\label{rep.odd}
\raisetag{10ex}
\psi^{(2n+1)}_{\lambda}(x_1)
\ket{k_1,\dots,k_n} &= \lambda
q^{-(k_1+\dots+k_n)}
\ket{k_1,\dots,k_n},   \\
\psi^{(2n+1)}_{\lambda}(x_1^*)
\ket{k_1,\dots,k_n} &=  \bar{\lambda}
q^{-(k_1+\dots+k_n)}
\ket{k_1,\dots,k_n},   \\
\psi^{(2n+1)}_{\lambda}(x_i)
\ket{k_1,\dots,k_n} &=  (1-
q^{-2k_{i-1}})^{1/2}
q^{-(\sum_{j=i}^{n}k_j)}
\ket{\dots,k_{i-1}-1,\dots },  \\
\psi^{(2n+1)}_{\lambda}(x_i^*) \ket{k_1,\dots,k_n} &=
(1-q^{-2(k_{i-1}+1)})^{1/2} q^{-(\sum_{j=i}^{n}k_j)}
\ket{\dots,k_{i-1}+1,\dots },
\end{split}
\end{equation}
for
$i=2,\dots,n+1$. With this values for the index $i$, $x_i$ lowers $k_{i-1}$,
whereas $x_i^*$ raises $k_{i-1}$. From  irreducibility the
vectors $\{ \ket{k_1,\dots,k_n} , k_i \geq 0 \}$ form an orthonormal
basis of the Hilbert space $\ch$.

As for the even case, the formul{\ae} \eqref{rep.odd} give bounded
operators only for $\abs{q}>1$; and as before, the representations for
$\abs{q}<1$ can be obtained by replacing all indices $k_i$ with
$-k_i-1$.

Again, as in the even case, we can verify  that
$\psi^{(2n+1)}_\lambda(x_1)$ is compact (indeed, trace class) and
that the  ideal
generated by $\psi^{(2n+1)}_\lambda(x_1)$ is  $\ck(\ch)$, in the \Cs-algebra
completion of the  image
$\psi^{(2n+1)}_\lambda(A(S^{2n+1}_q))$ . The
representations
$\psi^{(2n+1)}_\lambda$ can  be assembled into a single
representation by adjointable
operators on  a Hilbert $\cc(S^1)$-module. With this we can verify
that the ideal
generated by $x_1$ in $\cc(S^{2n+1}_q)$ is $\cc(S^1)\otimes\ck$ and this
verifies the
extension \eqref{ext.odd}.

\bigskip 
Summing up, we get a complete picture of the set of  irreducible representations of all these spheres
$S^{N}_q$; or equivalently, of the primitive spectrum of the \Cs-algebra $\cc(S^{N}_q)$ of continuous
functions on $S^{N}_q$.\\ For the odd spheres
$S^{2n+1}_q$, the set  of irreducible representations is indexed by the union of $n+1$ copies of
$S^1$.  These run from the representations $\psi^{(2n+1)}_\lambda$ of $S^{2n+1}_q$ given in
\eqref{rep.odd} down to the one dimensional representations $\psi^{(1)}_\lambda$ that factor through 
the undeformed $ S^1 $. \\
For the even spheres $S^{2n}_q$, the set of irreducible
representations is indexed by
the union of
$n$ copies of $S^1$ and $2$ points. The isolated points correspond to
the $2$
representations $\psi^{(2n)}_\pm$ specific to $S^{2n}_q$ and given in
\eqref{rep.even};
the circles  correspond to  representations $\psi^{(2m+1)}_\lambda$
coming from lower
odd  dimensional spheres, down to the undeformed $S^1$.

\section{$K$-homology and $K$-theory for Quantum Spheres}\label{se:fmqes}

We explicitly construct complete sets of generators for the $K$-theory (by
nontrivial self-adjoint idempotents and unitaries) and 
the $K$-homology (by nontrivial Fredholm modules) of the spheres
$S_q^{N-1}$. We also construct the corresponding Chern characters in cyclic
homology and cohomology and compute the  pairing of $K$-theory with
$K$-hom\-o\-logy.

We shall study generators of the  $K$-homology  and
$K$-theory of the spheres $S^{N-1}_q$.  The $K$-theory classes will be given by
means of self-adjoint idempotents (naturally associated with  the 
aforementioned
unipotents) and of unitaries in algebras of matrices over $A(S^{N-1}_q)$.
The $K$-homology classes will  be given as (homotopy classes of) suitable
$1$-summable Fredholm modules using the representations constructed previously.
\\
For odd spheres (i.e. for $N$ even) the odd $K$-homology generators are first given in terms of 
unbounded Fredholm modules. These are given by means of a natural unbounded operator $D$ which,
while failing to have compact resolvent, has bounded commutators with all
elements in the algebra $A(S^{2n-1}_q)$.

In fact, in order to compute the  pairing of $K$-theory with
$K$-homology, it is
more convenient to first compute the Chern characters and then use
the pairing
between
cyclic homology and cohomology \cite{co94}. Thus, together with the
generators
of
$K$-theory and $K$-homology we shall also construct the associated
Chern
characters in the cyclic homology $HC_*[A(S^{N-1}_q)]$ and cyclic
cohomology $HC^*[A(S^{N-1}_q)]$ respectively.

Needless to say, the pairing is integral (it comes from a
noncommutative index
theorem). The non-vanishing of the pairing will testify
to the
non-triviality  of the elements that we construct in both
$K$-homology $K$-theory.

\bigskip

It is worth recalling the $K$-theory
and homology of the
classical spheres.
\\
For an even dimensional sphere $S^{2n}$, the
groups are
\begin{gather*} K^0(S^{2n}) \cong \IZ^2 ,\qquad K^1(S^{2n}) =0 , \\
K_0(S^{2n})
\cong \IZ^2 ,\qquad K_1(S^{2n}) =0.
\end{gather*} One generator of the $K$-theory $[1]\in K^0(S^{2n})$ is
given by
the trivial $1$-dimensional bundle. The other generator of
$K^0(S^{2n})$ is the
left handed  spinor bundle. One $K$-homology generator $[\varepsilon]\in
K_0(S^{2n})$ is  ``trivial'' and is the push-forward of the generator
of $K_0(*)\cong\IZ$ by the inclusion
$\iota:*\into S^{2n}$ of a point (any point) into the sphere. The
other generator,
$[\mu]\in  K_0(S^{2n})$, is the $K$-orientation of
$S^{2n}$ given by
its
structure as a spin manifold \cite{co94}.

For an odd dimensional sphere, the groups are
\begin{gather*} K^0(S^{2n+1}) \cong \IZ ,\qquad K^1(S^{2n+1}) \cong \IZ
,
\\ K_0(S^{2n+1}) \cong \IZ ,\qquad K_1(S^{2n+1}) \cong \IZ .
\end{gather*} The generator $[1]\in K^0(S^{2n+1})$ is the
trivial
$1$-dimensional bundle. The generator of $K^1(S^{2n+1})$ is
a
nontrivial unitary matrix-valued function on $S^{2n+1}$; for instance, it may be
takes as the matrix \eqn{unit} in the limit $q=1$. The trivial
generator
$[\varepsilon]\in K_0(S^{2n+1})$  is again given
by the inclusion of a point. The generator $[\mu]\in
K_1(S^{2n+1})$ is the $K$-orientation of $S^{2n+1}$ given by its
structure as a spin manifold \cite{co94}.

There is a natural pairing between $K$-homology and $K$-theory. If we
pair
$[\varepsilon]$ with a vector bundle we get the rank of the vector
bundle, i.~e.\ the dimension of its fibers. If we pair $[\mu]$ with a
vector bundle it gives the ``degree'' of the bundle, a measure of its
nontriviality. Similarly, pairing with $[\mu]$ measures the
nontriviality of a unitary.

\bigskip
The $K$-theory and $K$-homology of the quantum Euclidean spheres are
isomorphic to that of the classical spheres; that is, for any $N$ and
$q$, one has that $K_*[\cc(S^{N-1}_q)] \cong K^*(S^{N-1})$ and $K^*[\cc(S^{N-1}_q)]
\cong K_*(S^{N-1})$.
In the case of $K$-theory, this was proven by Hong and Szyma{\'n}ski in
\cite{HonSzy01} using their construction of the \Cs-algebras as
Cuntz-Krieger algebras of graphs. The groups $K_0$ and $K_1$ were  given as the
cokernel and the kernel respectively, of a matrix canonically 
associated with
the graph.  The result for $K$-homology can be proven using the same techniques
\cite{Cun84,RaeSzy00}: the groups $K^0$ and $K^1$ are now given as 
the kernel and
the cokernel respectively, of the transposed matrix.
The $K$-theory and the $K$-homology for the particular case of $S^2_q$ (in fact
for all Podle{\'s} spheres $S^2_{qc}$) was worked out in \cite{MasNakWat91})
while for $S^3\cong SU_q(2)$ it was spelled out in \cite{MasNakWat90}.
\subsection{$K$-homology}
\label{se:khom}
Because the  $K$-homology of these deformed spheres is isomorphic to
the
$K$-hom\-o\-logy of the ordinary spheres, we need to construct  two
independent
generators. First consider the ``trivial'' generator of
$K^0[\cc(S^{N-1}_q)]$.  This can be constructed in a manner closely
analogous to the  undeformed case.

As we have already mentioned, the trivial generator of $K_0(S^{N-1})$
is  the
image of the  
generator of the $K$-homology of a point by the  functorial map
$K_*(\iota) : K_0(*) \to K_0(S^{N-1})$, where $\iota : *\into
S^{N-1}$
is the inclusion
of a point into the sphere. The  quantum Euclidean spheres do not
have as many points,
but they do have some.
We have seen that the relations among the various spheres always
include
a  homomorphism
$A(S_q^{N-1}) \to A(S^1)$. Equivalently, every $S^{N-1}_q$ has  a
circle $S^1$ as a
classical subspace;  thus for  every
$\lambda\in S^1$ there is a point, i.~e., the  homomorphism
$\psi^{(1)}_\lambda :
\cc(S^{N-1}_q) \to \IC$. \\
We can construct an element $[\varepsilon_\lambda] \in K^0[\cc(S_q^{N-1})]$ by pulling
back the generator  of $K^0(\IC)$ by $\psi^{(1)}_\lambda$. This construction factors
through $K_0(S^1)$. Because $S^1$ is path connected, the points of
$S^1$ all define homotopic (and hence $K$-homologous) Fredholm
modules. Thus there is a single $K$-homology class
$[\varepsilon_\lambda] \in K^0[\cc(S_q^{N-1})]$, independent of
$\lambda\in S^1$. 
\\
The canonical generator of $K^0(\IC)$ is  given by the following
Fredholm module: The
Hilbert space is $\IC$; the grading operator is $\gamma=1$; the
representation is the
obvious representation of $\IC$ on $\IC$; the Fredholm operator is
$0$. If we  pull this
back to $K^0[\cc(S^{N-1}_q)]$ using
$\psi^{(1)}_\lambda$, then  the Fredholm module
$\varepsilon_\lambda$ is given in the
same way  but with
$\psi^{(1)}_\lambda$ for the representation. 

Given this construction of $\varepsilon_\lambda$, it is straightforward to
compute its Chern
character
$\ChC^*(\varepsilon_\lambda) \in HC^*[A(S^{N-1}_q)]$: It is the pull  back of
the Chern character
of the canonical generator of $K^0(\IC)$.  An element of the cyclic
cohomology
$HC^0$ is a  trace. The degree $0$ part of the
Chern character of the canonical  generator of $K^0(\IC)$ is given by
the identity map
$\IC\to\IC$,  which is trivially a trace. Pulling this back we find
$\ChC^0(\varepsilon_\lambda) = \psi^{(1)}_\lambda : A(S^{N-1}_q)
\to \IC$ which is also  a trace because it is a homomorphism to a
commutative algebra. These are distinct elements of $HC^0[A(S^{N-1}_q)]$ for  different values of
$\lambda$. However, because the Fredholm modules $\varepsilon_\lambda$ all lie in the same
$K$-homology class, their Chern characters are all equivalent in {\it periodic} cyclic cohomology
defined in \eqn{pecy}. Indeed, applying the periodicity operator \eqn{peop} once one gets that the
cohomology classes $S(\psi^{(1)}_\lambda) \in HC^2[A(S^{N-1}_q)]$ are all the same. For the computation
of the pairing between $K$-theory and
$K$-homology, any trace determining the same periodic cyclic cohomology class can be used. The most
symmetric choice of trace is given by  averaging $\psi^{(1)}_\lambda$ over $\lambda\in S^1\subset
\IC$:
\[
\tau^0(a) := \oint_{S^1} \psi^{(1)}_\lambda(a) \frac{d\lambda}{2\pi i
\lambda} .
\]
The result is normalized, $\tau^0(1)=1$, and vanishes
on all the generators. The higher degree parts of $\ChC^*(\varepsilon_\lambda)$
depend only on the $K$-homology class $[\varepsilon_\lambda]$ and can be
constructed from $\tau^0$ by the periodicity
operator \eqn{peop}.

\subsection{Fredholm Modules for even spheres}
We will now construct an
element $[\mu_{\mathrm{ev}}] \in K^0[\cc(S^{2n}_q)]$ by giving a
suitable even
Fredholm module
$\mu := (\ch,F,\gamma)$. \\
Identify the Hilbert spaces for
the representations $\psi^{(2n)}_\pm$ given in \eqref{rep.even} by
identifying their bases, and
call this $\ch$. The representation for
the Fredholm module
is
\[
\psi:=\psi^{(2n)}_+\oplus\psi^{(2n)}_-
\] acting on
$\ch\oplus\ch$. The grading operator and the Fredholm
operator are respectively,
\[
\gamma = \begin{pmatrix} 1&0
\\0&-1
\end{pmatrix}
,\quad
  F = \begin{pmatrix} 0&1\\1&0
\end{pmatrix}.
\]
It is obvious that $F$ is odd (since it anticommutes with $\gamma$)
and Fredholm (since
it is invertible).
The remaining property to check is that for any $a\in A(S^{2n}_q)$,
the commutator $[F,\psi(a)]_-$ is compact. Indeed,
\[
[F,\psi(a)]_- =
\begin{pmatrix} 0 & - \psi^{(2n)}_+(a) + \psi^{(2n)}_-(a) \\
\psi^{(2n)}_+(a) - \psi^{(2n)}_-(a) & 0
\end{pmatrix} .
\]
However, $\psi^{(2n)}_+(a) - \psi^{(2n)}_-(a) =
\psi^{(2n)}_+[a -
\sigma(a)]$ and $a-\sigma(a)$ is always proportional to a power of
$x_0$. Thus   this
is not only compact, it is trace class. This also shows that we have
(at least) a $1$-summable Fredholm module. This is in contrast to the fact that the analogous 
element of $K_0(S^{2n})$ for the undeformed sphere is given by a $2n$-summable Fredholm module.
\\
The corresponding Chern character \cite{co94} $\ChC^*(\mu_{\mathrm{ev}})$ has  a
component
in degree $0$, $\ChC^{0}(\mu_{\mathrm{ev}})\in HC^0[A(S^{2n}_q)]$. 
From the general construction
\eqn{chfreven}, the element
$\ChC^0(\mu_{\mathrm{ev}})$ is the trace
\begin{equation}\label{tau1}
\tau^1(a) := \tfrac12 \Tr\left(\gamma F ([F,\psi(a)] \right) = \Tr\left[\psi^{(2n)}_+(a) -
\psi^{(2n)}_-(a)\right].
\end{equation}
As we have mentioned, $\psi^{(2n)}_+(a) - \psi^{(2n)}_-(a) = \psi^{(2n)}_+[a - \sigma(a)]$ is trace class
since $a-\sigma(a)$ is always proportional to a power of
$x_0$. The higher degree parts of $\ChC^*(\mu_{\mathrm{ev}})$
can be obtained via the periodicity operator \eqn{peop}.

For  $S^{2}_q$ our Fredholm module coincides with the one 
constructed in
\cite{MasNakWat91}.

\subsection{Fredholm Modules for odd spheres}\label{se:khodd}

The element $[\mu_{\mathrm{odd}}] \in
K^1[\cc(S^{2n+1}_q)]$ is most
easily given by an unbounded Fredholm
module. The corresponding unbounded operator $D$
which, while failing to have compact resolvent, has bounded 
commutators with all elements in the algebra $A(S^{N-1}_q)$.\\
Let the representation $\psi$ be
the direct integral (over $\lambda\in S^1$) of the
representations
$\psi^{(2n+1)}_\lambda$ given in \eqref{rep.odd}. The operator is the
unbounded `Dirac' operator
$D:= \lambda^{-1}\frac{d}{d\lambda}$.
\\
{}From \eqref{rep.odd}, we see that the representative of $x_1$
is proportional to $\lambda$ and as a consequence,
\begin{subequations}
\begin{equation}
[D,\psi(x_1)]_- = \psi(x_1)
\end{equation}
whereas for $i>1$ the representative of $x_i$ does not involve
$\lambda$ and therefore
\begin{equation}
[D,\psi(x_i)]_- = 0, \quad i>0 .
\end{equation}
\end{subequations}
Since $a\mapsto [D,\psi(a)]_-$ is a derivation, this shows that
$[D,\psi(a)]_-$ is bounded for any $a \in A(S^{2n+1}_q)$. Note
however that for $n>0$ (i.~e., except for $S^1$) all eigenvalues of
$D$ have infinite degeneracy and therefore $D$ does not have compact
resolvent.

This triple can be converted in to a bounded  Fredholm module by
applying a cutoff
function to $D$. A convenient  choice is $F=\chi(D)$ where
\[
\chi(m) :=
\begin{cases} 1 &: m>0\\ -1  &: m\leq 0 .
\end{cases}
\] To be more explicit, use a Fourier  series basis for the Hilbert
space,
\[
\ket{k_0,k_1,\dots,k_n} :=
\lambda^{k_0}\ket{k_1,\dots,k_n} ,
\]
in which the  representation is given  by,
\begin{align*}
\psi(x_1)
\ket{k_0,\dots,k_n} &=
q^{-(k_1+\dots+k_n)}
\ket{k_0+1,\dots,k_n}
,  \\
\psi(x_1^*)
\ket{k_0,\dots,k_n} &=
q^{-(k_1+\dots+k_n)}
\ket{k_0-1,\dots,k_n},  \\
\psi(x_i)
\ket{k_0,\dots,k_n} &=  (1-q^{-2k_{i-1}})^{1/2} q^{-(k_i+\dots+k_n)}
\ket{ \dots,k_{i-1}-1,\dots }, \\
\psi(x_i^*) \ket{k_0,\dots,k_n} &= (1-q^{-2(k_{i-1}+1)})^{1/2}
q^{-(\sum_{j=i}^{n}k_j)}
\ket{ \dots,k_{i-1}+1,\dots },
\end{align*}
for $i = 1, \dots, n$. The Fredholm operator is then given by
\[
F\ket{k_0,\dots,k_n} =
\chi(k_0) \ket{k_0,\dots,k_n}
.\]
The only condition to  check is that the commutator $[F,\psi(a)]_-$
is compact for any
$a\in  \cc(S^{2n+1}_q)$. Since $a \mapsto [F,\psi(a)]_-$ is a
derivation, it  is
sufficient to check this on generators. One  finds
\begin{align}
\label{F.com} 
& [F,\psi(x_i)]_-=0~, ~~ i>1~,  \nn
& [F,\psi(x_1)]_-
\ket{k_0,\dots,k_n} =
\begin{cases} 2q^{-(k_1+\dots+k_n)}
\ket{1,k_1,\dots,k_n} &:
k_0=0\\ 0 &: k_0\neq
0,
\end{cases}
\end{align} which is indeed compact, and in
fact trace class.
\\
Thus, this is a $1$-summable Fredholm module. Again this is in contrast to the fact that the 
analogous element of
$K_1(S^{2n+1})$ for the undeformed sphere  is given by a $(2n+1)$-summable Fredholm module.
\\
Its Chern character \cite{co94} begins with
$\ChC^{\frac12}(\mu_{\mathrm{odd}})\in HC^1[A(S^{2n+1}_q)]$. From the general construction
\eqn{chfrodd}, the element $\ChC^{\frac12}(\mu_{\mathrm{odd}})$ is
given by the cyclic $1$-cocycle $\varphi$ defined
by
\begin{equation}
\label{phi.def}
\varphi(a,b) := \tfrac12
\Tr\left(\psi(a)[F,\psi(b)]_-\right)
.\end{equation}
One checks directly cyclicity, i.e. $\varphi(a,b)=-\varphi(b,a)$, and closure under $b$, i.e.
$\varphi(ab,c)-\varphi(a,bc)+\varphi(ca,b)=b\varphi(a,b,c)$.\\
The higher degree parts of $\ChC^*(\mu_{\mathrm{odd}})$
can be obtained via the periodicity operator \eqn{peop}.

For $S^3_q \cong \SU_q(2)$ our Fredholm module coincides with the
one  in
\cite{MasNakWat90}.

\subsection{Singular integrals}
We could interpret the classes 
$[\mu_{\mathrm{ev}}] \in K^0[\cc(S^{2n}_q)]$ and $[\mu_{\mathrm{odd}}] \in K^1[\cc(S^{2n+1}_q)]$ as
giving `singular' integrals over the corresponding quantum spheres. With the associated Chern
characters given in \eqn{tau1} and \eqn{phi.def} respectively, and from the general
expression  \eqn{chfr}, these integrals  are given by,
\begin{align}
 \int_{S^{2n}_q} a &= \tau^1(a) ~, 
~~~ \forall a\in \ca(S^{2n}_q)~, \\
 \int_{S^{2n+1}_q} a ~d b &= \phi(a,b) ~, 
~~~ \forall a,b\in \ca(S^{2n+1}_q)~.
\end{align} 
As a way of illustration, let us compute them on
generators. We indicate with $\delta_{ij}$ the usual Kronecker delta which is
equal to $1$ if
$i=j$ and $0$ otherwise.\\ Firstly, for even spheres we find
\begin{equation}
\int_{S^{2n}_q} x_i  = \tau^1(x_i) = 
\Tr\left[\psi^{(2n)}_+\left(x_i - \sigma(x_i)\right)\right] = 2 \Tr\left[\psi^{(2n)}_+(x_i)\right] 
\delta_{i0} .
\end{equation}
Thus, we need to compute
\begin{align*}
\Tr[\psi^{(2n)}_+(x_0)] &=
\sum_{k_0=0}^\infty \dots \sum_{k_{n-1}=0}^\infty
q^{-(k_0+\dots+k_{n-1})} =
\left(\sum_{k=0}^\infty q^{-k}\right)^n \\ &= (1- q^{-1})^{-n}
,\end{align*} 
and in turn, 
\begin{equation}
\int_{S^{2n}_q} x_i = \frac{2}{(1-q^{-1})^n}
~\delta_{i0}~.  
\end{equation}
Similarly, for odd spheres we find
\begin{equation}
\int_{S^{2n+1}_q} x_i ~dx^*_j = \phi(x_i,x_i^*) = 
\frac{1}{2} \Tr\left(\psi(x_1^*)[F,\psi(x_1)]_-\right) \delta_{i1} \delta_{j1}.
\end{equation}
We have already computed
$[F,\psi(x_1)]_-$ in eq.~\eqref{F.com}. From  that, we get
\[
\psi(x_1^*)[F,\psi(x_1)]_- \ket{k_0,\dots,k_n}  = \begin{cases} 2
q^{-2(k_1+\dots+k_n)} \ket{0,k_1,\dots,k_n} &: k_0=0\\ 0 &:
k_0\neq 0 .
\end{cases}
\]
Thus,
\begin{align*}
\Tr\left(\psi(x_1^*)[F,\psi(x_1)]_-\right)  &= \sum_{k_1=0}^\infty
\dots
\sum_{k_n=0}^\infty 2 q^{-2(k_1+\dots+ k_n)} = 2
\left(\sum_{k=0}^\infty q^{-2k}\right)^n \\ &= 2 (1- q^{-2})^{-n}
,\end{align*}
and in turn, 
\begin{equation}
\int_{S^{2n+1}_q} x_i ~dx^*_j = \frac{1}{(1-q^{-2})^n} ~\delta_{i1}\delta_{j1}~.  
\end{equation}

\subsection{$K$-theory for even  spheres}
For $S^{2n}_q$ we construct
two classes in the $K$-theory group
$K_0[\cc(S^{2n}_q)] \cong \IZ^2$. \\
The first class is trivial. The element $[1]\in
K_0[\cc(S^{2n}_q)]$ is
the equivalence class of $1\in \cc(S^{2n}_q)$
which is of course an
idempotent. In order to compute the pairing
with $K$-homology, we need
the degree
$0$ part of its Chern
character,
$\ChC_0[1]$, which is represented by the cyclic cycle $1$. \\
The second, nontrivial, class was presented in \cite{LanMad01}.
It is  given by an idempotent
$e_{(2n)}$ constructed from the unipotent \eqref{uni.even} as
\begin{equation}
\label{proj} e_{(2n)} =
\tfrac12(\II+u_{(2n)}) 
\end{equation}
(again, for the sphere $S^{2}_q$ the idempotent
\eqref{proj}  was already  
in \cite{MasNakWat91}). Its degree $0$ Chern character, $\ChC_{0}(e_{(2n)})\in
HC_0[A(S^{2n}_q)]$, is
\begin{align}
\ChC_0(e_{(2n)}) &= \tr(e_{(2n)} - \tfrac{1}{2} \II_{2n}) = \tfrac12 \tr (u_{(2n)}) \nn
 &= \tfrac12(q^{-1}-1)^n x_0 ,
\end{align}
since the recursive definition \eqref{uni.even} of the unipotent
$u_{(2n)}$ shows that,
\[
\tr(u_{(2n)}) = (q^{-1}-1) \tr(u_{(2n-2)}) = (q^{-1}-1)^n x_0 .
\]

Now, we can pair these classes with the two $K$-homology elements
which we
constructed in Section~\ref{se:khom}. First,
\[
\langle \varepsilon_\lambda, [1]\rangle := \tau^0(1) =
1,
\]
which is hardly surprising. Second, the ``rank'' of  the idempotent
$e_{(2n)}$ is
\[
\langle \varepsilon_\lambda, e_{(2n)} \rangle
:= \tau^0(\tr(e_{(2n)}) = 2^{n-1}
.\]
Also not surprising is the ``degree'' of $[1]$,
\[
\langle \mu_{\mathrm{ev}}, [1]\rangle := \tau^1(1) =
\Tr\bigl[\psi^{(2n)}_+(1)-\psi^{(2n)}_-(1)\bigr] =
\Tr(1-1) = 0
.\]
The more complicated pairing is,
\begin{align*}
\langle \mu_{\mathrm{ev}}, e_{(2n)}\rangle  &:= \tau^1(\ChC_0
e_{(2n)}) \\ &\:=
\Tr\circ\psi^{(2n)}_+\circ(1-\sigma)\left(2^{n-1} + \tfrac12
[q^{-1}-1]^nx_0\right) \\ &\:= (q^{-1}-1)^n \Tr[\psi^{(2n)}_+(x_0)] 
 = (q^{-1}-1)^n (1-q^{-1})^{-n} \\
&\:= (-1)^n
.\end{align*} 
The fact that the matrix of pairings,
\[
\begin{array}{c|cc}
& [1] & [e_{(2n)}] \\
\hline
[\varepsilon_\lambda] & 1 & 2^{n-1} \\
{}[\mu_{\mathrm{ev}}] & 0 & (-1)^n
\end{array}
\]
is invertible over the integers proves that the classes
$[1],[e_{(2n)}]$, both elements of the group $K_0[\cc(S^{2n}_q)] \cong \IZ^2$, and the classes
$[\varepsilon_\lambda],[\mu_{\mathrm{ev}}]$, both in $K^0[\cc(S^{2n}_q)] \cong \IZ^2$, 
are nonzero and that no one of them  may be a multiple of another class; thus they  are generators of
the respective groups.

Classically, the ``degree'' of the left-handed spinor
bundle is
$-1$. So, the $K$-hom\-o\-logy class which correctly generalizes the classical
$K$-orientation class $[\mu]\in K_0(S^{2n})$ is actually
$(-1)^{n+1}[\mu_{\mathrm{ev}}]$.

\subsection{$K$-theory for odd spheres}
Again, define
$[1]\in
K_0[\cc(S^{2n+1}_q)]$ as the equivalence class of
$1\in
\cc(S^{2n+1}_q)$. The pairing with our element $[\varepsilon_\lambda]
\in
K^0[\cc(S^{2n+1}_q)]$ is again,
\[
\langle \varepsilon_\lambda, [1]\rangle
:=\tau^0(1) = 1
.\]
There is no other independent
generator in $K_0[\cc(S^{2n+1}_q)] \cong \IZ$.

Instead, $K_1[\cc(S^{2n+1}_q)] \cong \IZ$ is nonzero. So we need to
construct a
generator there. An odd $K$-theory element is an equivalence class of
unitary
matrices over the algebra. We can construct an appropriate sequence
of unitary matrices recursively, just as we constructed the
unipotents and
idempotents. \\
Let $V_{(2n+1)} \in \Mat_{2^n}(A(S^{2n+1}_q))$ be defined recursively by
\begin{equation} \label{unit}
V_{(2n+1)} =
\begin{pmatrix} x_{n+1}  & q^{-1}
V_{(2n-1)} \\ -V_{(2n-1)}^* &
x_{n+1}^*
\end{pmatrix}
,
\end{equation}
with
$V_{(1)} = x_1$. By using the defining relations
\eqref{core} one
directly proves that it is unitary: 
\begin{equation}
V_{(2n+1)}V_{(2n+1)}^* =V_{(2n+1)}^* V_{(2n+1)} = 1~.
\end{equation}

In order to pair our
$K$-homology element $[\mu_{\mathrm{odd}}]\in
K^1[\cc(S^{2n+1}_q)]$ with the unitary $V_{(2n+1)}$, we
need the lower degree part $\ChC_{\frac12}(V_{(2n+1)}) \in
HC_1[A(S^{2n+1}_q)]$ of its Chern character. It is
given by the cyclic cycle,
\begin{align}
\ChC_{\frac12}(V_{(2n+1)}) &:=
\tfrac12\tr\left(V_{(2n+1)} \otimes V_{(2n+1)}^* -
V_{(2n+1)}^*
\otimes V_{(2n+1)} \right) \nn  & \:=
\tfrac12(q^{-2}-1)^n  (x_1
\otimes x_1^* - x_1^*\otimes x_1) .
\end{align}
Now, compute the pairing,
\begin{align*}
\langle \mu_{\mathrm{odd}} , V_{(2n+1)}\rangle  &:= \langle
\varphi,
\ChC_{\frac12}(V_{(2n+1)})\rangle \\ &\:= - (q^{-2}-1)^n
\varphi(x_1^*,x_1) \\ &\:= -
\tfrac12(q^{-2}-1)^n
\Tr\left(\psi(x_1^*)[F,\psi(x_1)]_-\right) 
 \\ &\:= -
\tfrac12(q^{-2}-1)^n 2 (1- q^{-2})^{-n} \\
&\:= (-1)^{n+1}
.\end{align*}
This proves that  $[V_{(2n+1)}]\in K_1[\cc(S^{2n+1}_q)]$ and
$[\mu_{\mathrm{odd}}]\in K^1[\cc(S^{2n+1}_q)]$ are nonzero and that
neither may be a multiple of another class. Thus $[V_{(2n+1)}]$ and
$[\mu_{\mathrm{odd}}]$ are indeed generators of these groups.

\subsection*{Acknowledgments.}
These lecture notes are based on work done with Alain Connes, Eli Hawkins, John Madore and Joe
V\'arilly; I am most grateful to them.   I thank Satoshi Watamura, Yoshiaki Maeda and Ursula
Carow-Watamura for their kind invitation to Sendai and for the fantastic hospitality there. I also thank
all participants of the conference for the great time we had together and Tetsuya Masuda for his help and
hospitality in Sendai and Tokyo. Finally, I thank Ludwik D\c{a}browski, Eli Hawkins, Fedele
Lizzi and Denis Perrot for very useful suggestions which improved the compuscript.

\vspace{1 cm}
\providecommand{\href}[2]{#2}\begingroup\raggedright
\endgroup

\end{document}